\documentclass[12pt]{article}
\usepackage{amsmath}
\usepackage{amsfonts,amssymb}
\usepackage{epsfig}
\usepackage{psfrag}
\usepackage{wrapfig}
\usepackage{graphicx}
\newcommand{\Z}{{\mathbb Z}}
\newcommand{\mb}{\mathbf}
\newcommand{\B}{\bigskip}
\newcommand{\m}{\medskip}
\newcommand{\proof}{\noindent{\em Proof: }}

\newcommand{\qed}{\hspace{\fill}$\square$}
\newtheorem{theorem}{Theorem}
\newtheorem{cor}[theorem]{Corollary}
\newtheorem{proposition}[theorem]{Proposition}
\newtheorem{lemma}[theorem]{Lemma}
\newtheorem{alg}[theorem]{Procedure}
\newtheorem{remark}[theorem]{Remark}
\newtheorem{example}[theorem]{Example}

\def\sqr#1#2{{\vcenter{\vbox{\hrule height.#2pt
 \hbox{\vrule width.#2pt height#1pt \kern#1pt
 \vrule width.#2pt}
 \hrule height.#2pt}}}}
\def\square{\mathchoice\sqr68\sqr68\sqr{2.1}3\sqr{1.5}3}
\date{\today}

\title{
Tree Orbits under Permutation Group Action: Algorithm, Enumeration and
 Application to Viral Assembly
}

\author{
Mikl\'os B\'ona\footnote{supported in part by 
\uppercase{NSF} grant DMS0714912
Department of Mathematics, University of Florida, Gainesville, FL 32611.
Email: bona@math.ufl.edu}\\
Meera Sitharam\footnote{corresponding author; 
supported in part by \uppercase{NSF} 
grant DMSO714912;
CISE dept., University of Florida, Gainesville, FL 32611.
Email: sitharam@cise.ufl.edu}\\
Andrew Vince\footnote{Department of  Mathematics, University of Florida,
 Gainesville, FL 32611. Email: avince@math.ufl.edu.
}}

\begin{document}
\maketitle
\begin{abstract}

This paper uses combinatorics and group theory to answer questions
about the assembly of icosahedral viral shells.  Although the
geometric structure of the capsid (shell) is fairly well understood in
 terms of 
its constituent subunits, 
the assembly process is not. For the purpose of this paper, 
the capsid is modeled by a polyhedron whose
facets represent the monomers.  The
assembly process is modeled by a rooted tree, the leaves representing
the facets of the polyhedron, the root representing the assembled
polyhedron, and the internal vertices representing intermediate stages of
assembly (subsets of facets).  Besides its virological motivation,
the enumeration of orbits of trees under the action of a finite group is
of independent mathematical interest.  If $G$ is a finite group acting
on a finite set $X$, then there is a natural induced action of $G$ on the set
$\mathcal{T}_X$ of trees whose leaves are bijectively labeled by the
elements of $X$.  If $G$ acts simply on $X$, then $|X| := |X_n| = n
\cdot |G|$, where $n$ is the number of $G$-orbits in $X$.  The basic
combinatorial results in this paper are 
(1) a formula for the number of orbits of each
size in the action of $G$ on $\mathcal{T}_{X_n}$, for every $n$, and 
(2) a simple  algorithm to find the
stabilizer  of a tree $\tau \in \mathcal{T} _X$
 in $G$ that runs in linear time and
does not need memory in addition to its input tree. 

\medskip\noindent
2000 Mathematics Subject Classification: Primary 05C05, 05A15,
20B25, 92C50. \vskip 1mm

Key words: tree enumeration, generating function, group action,
viral capsid assembly.

\end{abstract}


\section{Introduction}
\label{intro}

Viral shells, called {\it capsids}, encapsulate and protect the
fragile nucleic acid genome from physical, chemical, and enzymatic
damage.  Francis Crick and James Watson (1956) were the first to
suggest that viral shells are composed of numerous identical protein
subunits called {\it monomers}.  For many viruses, these monomers are
arranged in either a helical or an icosahedral structure. We are
interested in those shells that possess icosahedral symmetry.

Icosahedral viral shells can be classified based on their polyhedral
structure, facets corresponding to the monomers. The
classical ``quasi-equivalence theory'' of Caspar and Klug \cite{KC}
explains the structure of the polyhedral shell in the case where the
monomers have very similar neighborhoods.
  According to the theory, the number of
facets in the polyhedron is $60T$, where the $T$-{\it number} is of
the form $h^2 + hk + k^2$. Here $h$ and $k$ are non-negative
integers. The icosahedral group acts simply on the set of facets
of the polyhedron (monomers of the shell).

\begin{figure}
\centerline{
\psfig{file= 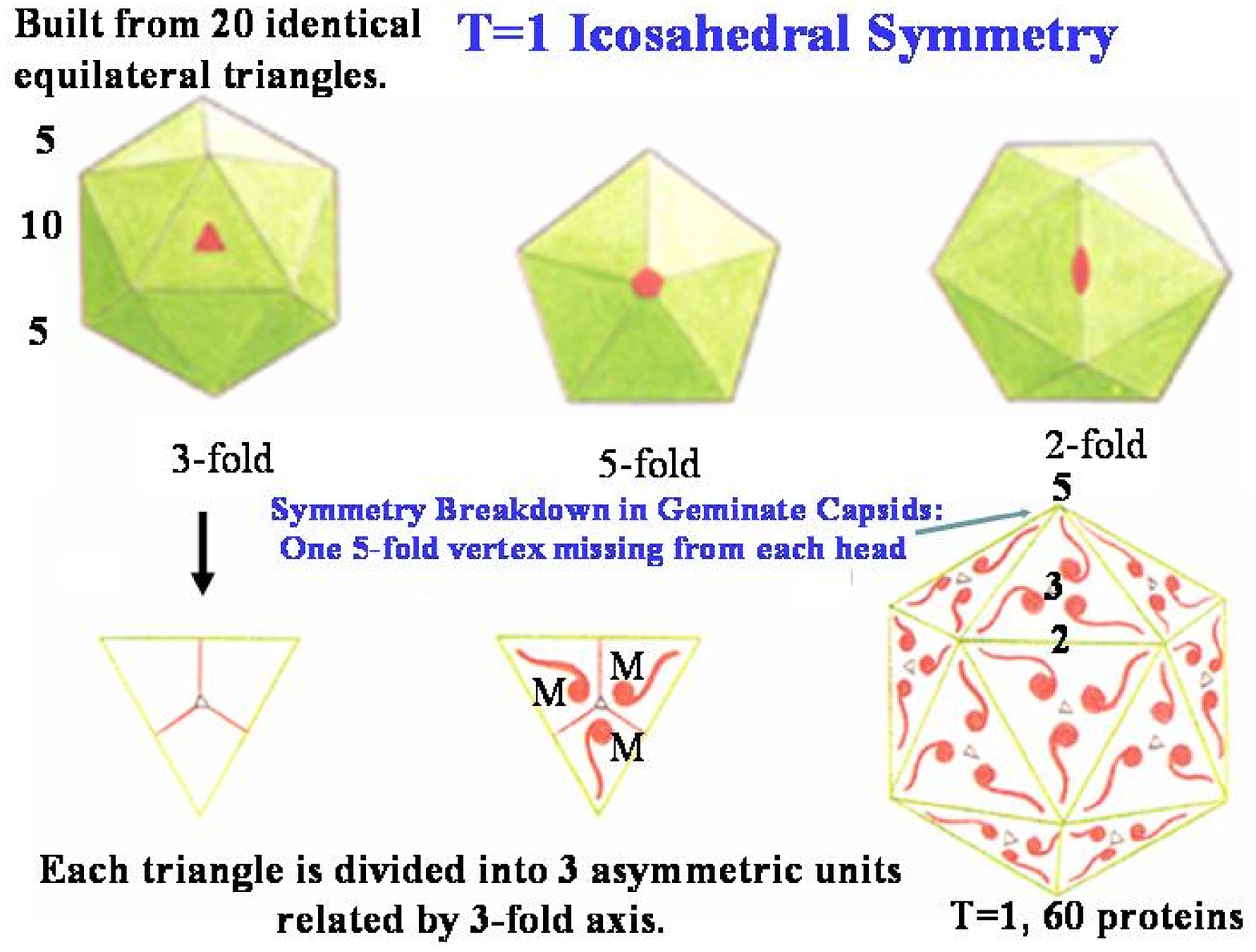, width = 7.8cm},
        \psfig{file=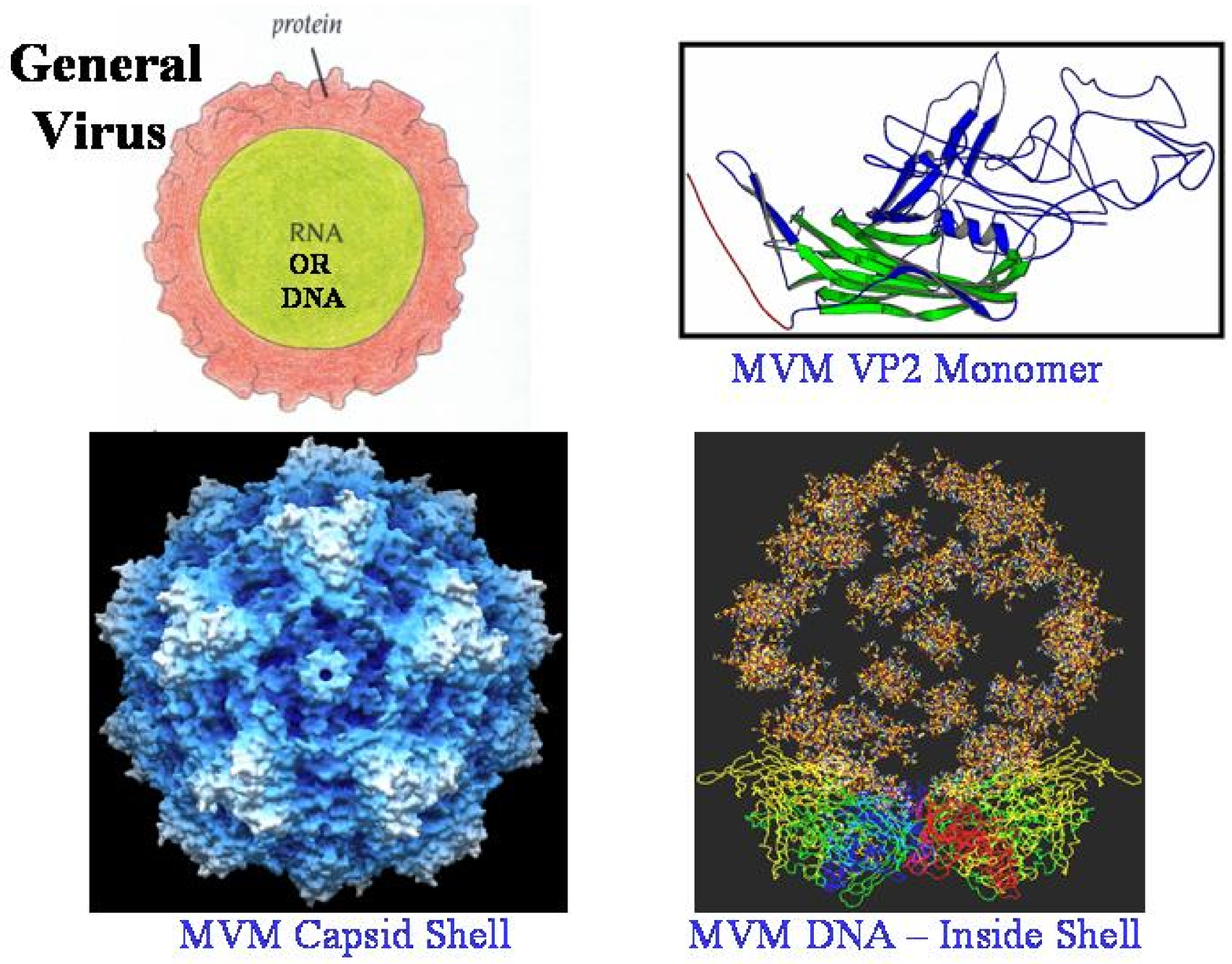, width=7.8cm}}
\caption{(Left) Basic Viral Structure. (Right)
   Minute virus of Mice X-ray and monomer structure courtesy 
\cite{McKenna-lab}}.
\label{virus-basic1-mvm1}
\end{figure}

Although for many virus families the structure is fairly well
understood and substantiated by crystallographic images, the viral
assembly process - just like many other spontaneous macromolecular
assembly processes - is not well-understood, even for $T=1$ viral
shells.  In many cases, the capsid self-assembles spontaneously,
rapidly and quite accurately in the host cell, with or without
enclosing the internal genomic material, and without the use of
chaperone, scaffolding or other helper proteins.  This is the type of
assembly that we consider here. See Figure \ref{virus-basic1-mvm1} for
basic icosahedral structure and X-ray structure of a T=1 virus.

\begin{figure}
\centerline{
\epsfig{file= 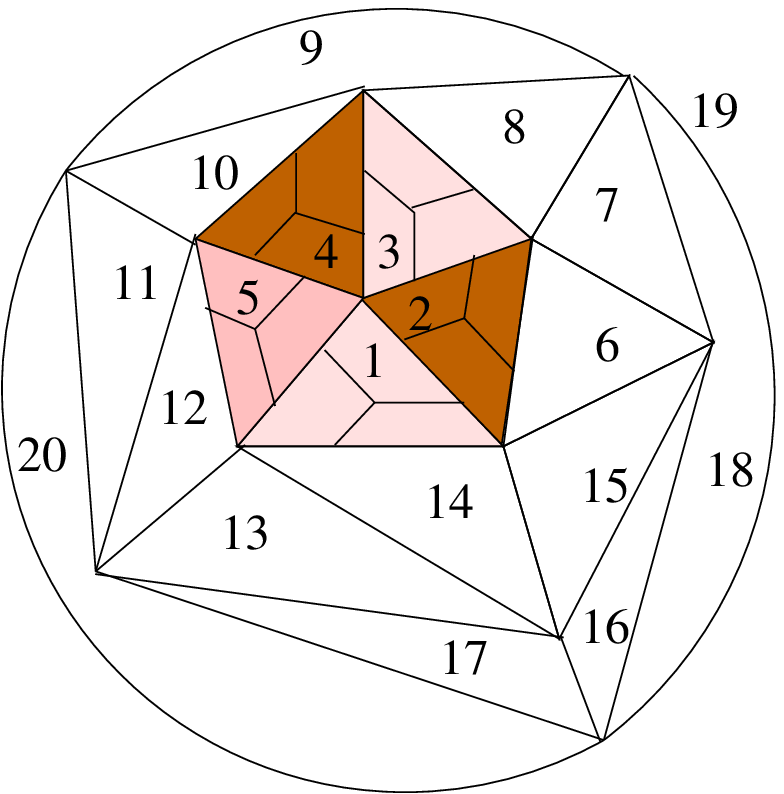, width = 7.3cm}
\epsfig{file= 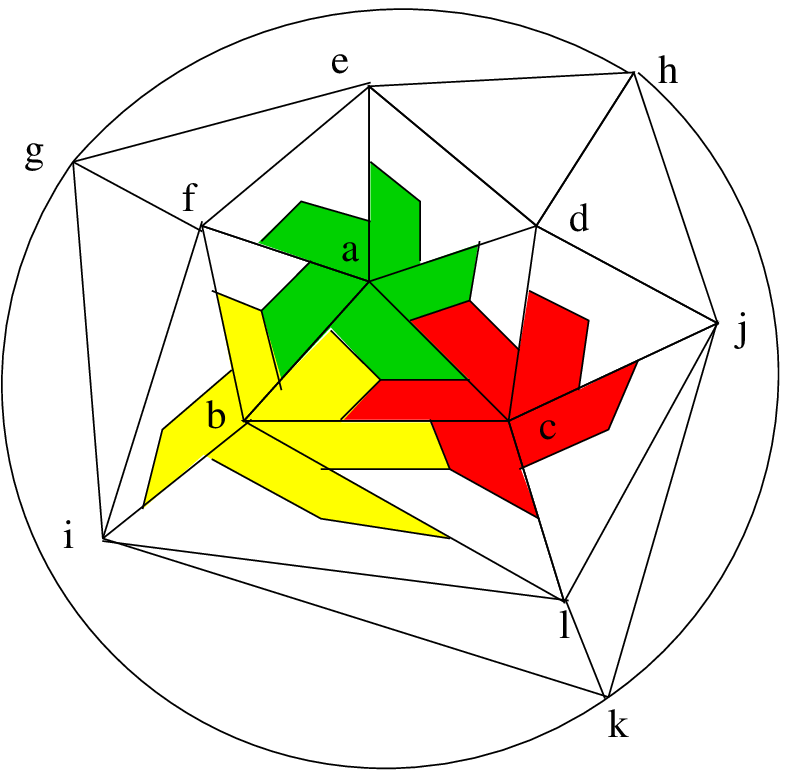, width = 7.3cm}
}
\caption{(Left) Face numbers: pentamer of trimers in a $T=1$ polyhedron.
(Right)
Vertex numbers: trimers of  pentamers in a $T=1$ polyhedron.
}
\label{virus-face-virus-vertex}
\end{figure}

\begin{figure}
\centerline{
\epsfig{file= 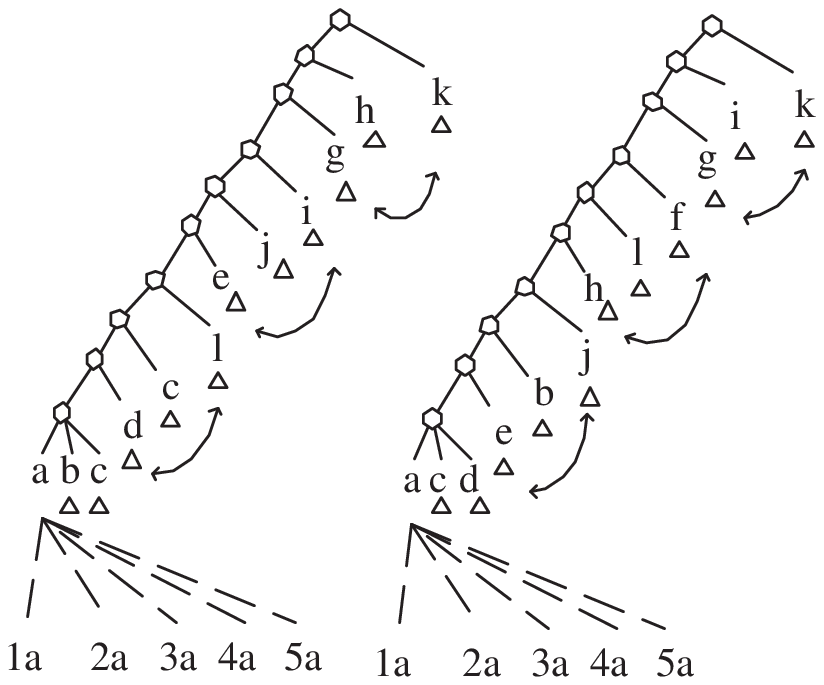, width = 8cm}
\epsfig{file= 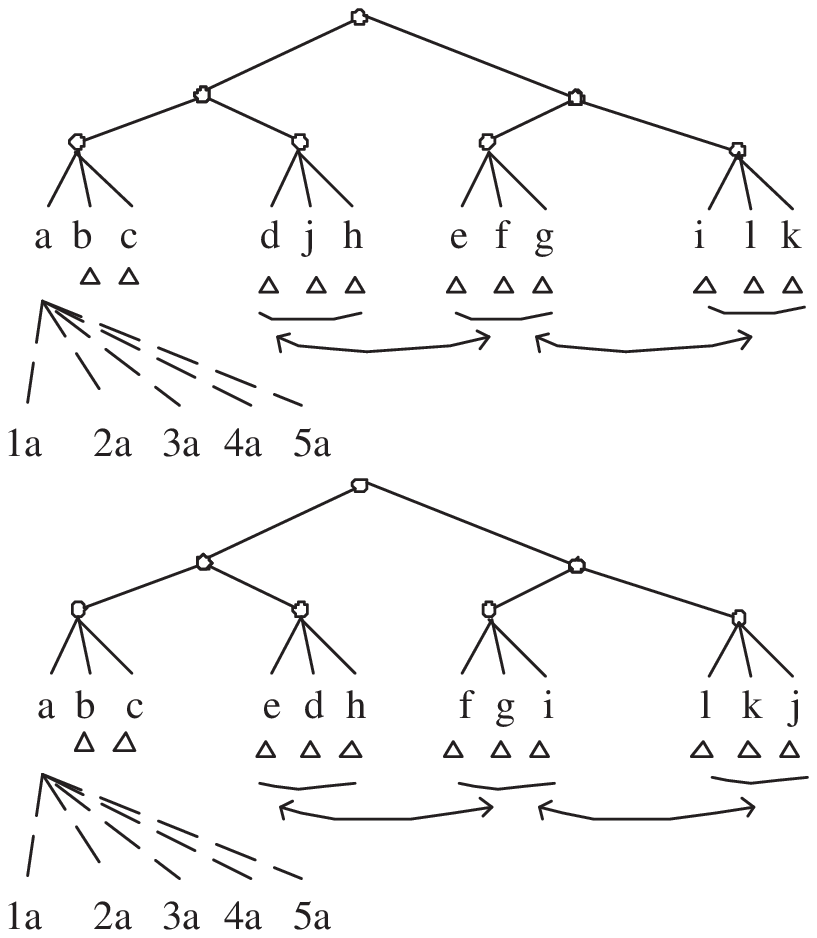, width = 7.5cm}
}

\caption{$T=1$ valid assembly trees based on pentameric subassemblies
and nomenclature of Figure \ref{virus-face-virus-vertex}, triangles
at bottom represent vertices with five leaves as children;
 arrows represent the action of 
the icosahedral group on trees ; only the long horizontal arrows
in the two figures on right fix the corresponding trees.}
\label{pathways1}
\end{figure}

\begin{figure}
\centerline{
\epsfig{file= 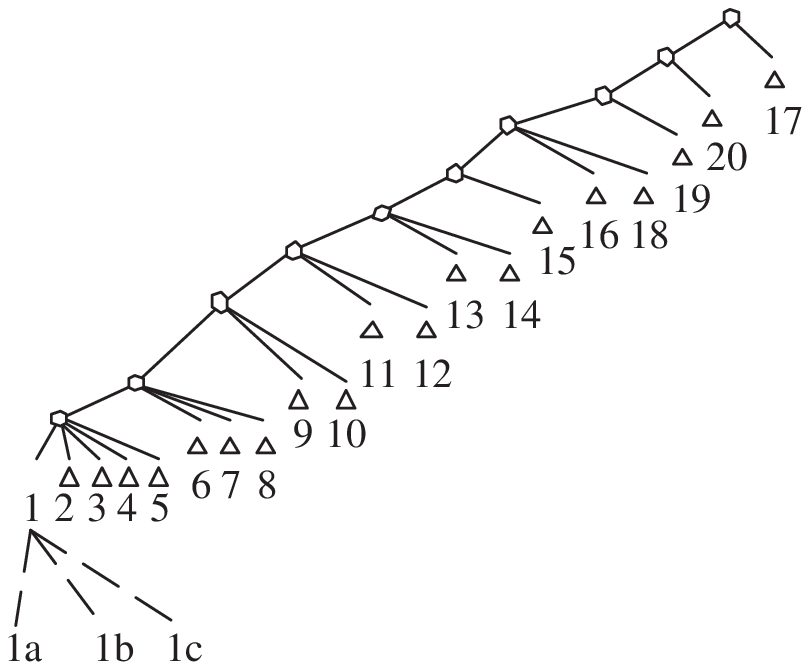, width = 6cm}
\epsfig{file= 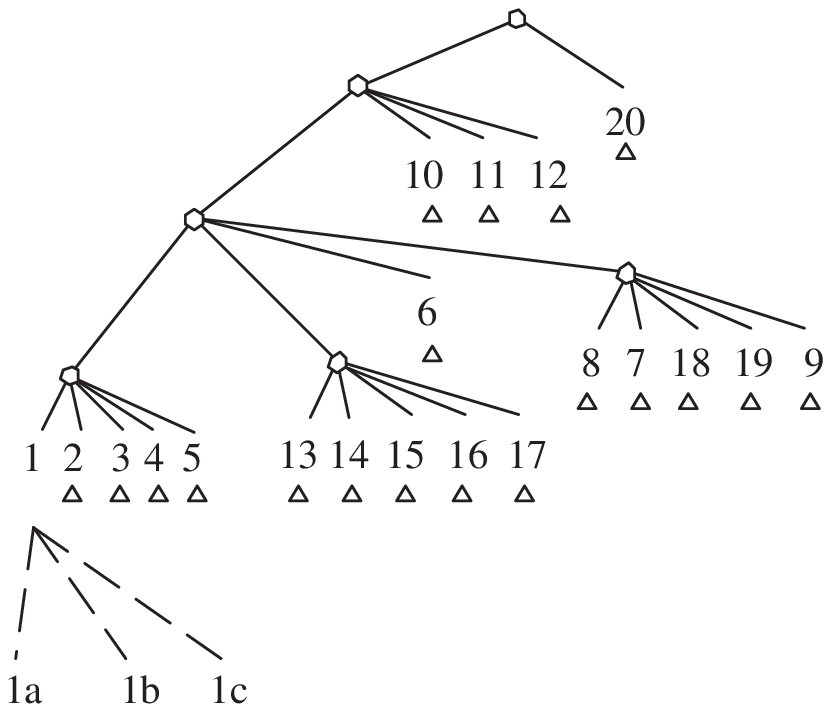, width = 6cm}
}
\caption{$T=1$ valid assembly trees based on trimeric subassemblies, triangles
at bottom represent trimers.}
\label{pathways2}
\end{figure}

Many mathematical models of viral shell assembly have been proposed
and studied including
\cite{berger2, bib:Berger,johnson-assembly2, marzec-day,
johnson-assembly,reddy, M1,zlotnick2,zlotnick1}.  Here we
use the {\it GT (geometry and tensegrity)} model of \cite{bib:AgSi06}. 
 In the GT model,
information about the construction (or decomposition) of the viral
shell is represented by an {\it assembly tree}.  The vertices of the
tree represent subassemblies that do not disintegrate during the
course of the assembly process.  In an assembly
tree, these subassemblies are partially ordered by containment, with
the root representing the complete assembled structure, and the leaves
representing the monomers.  That is, we only consider trees that
represent successful assemblies.  See Figure
\ref{virus-face-virus-vertex} for the nomenclature of a $T=1$ polyhedron,
and Figures \ref{pathways1},
\ref{pathways2} for examples of assembly trees.
Besides being intuitive and analyzable, it was shown in
\cite{bib:AgSi06} that the GT model's rough predictions fit
experimental and biophysical observations of known $T=1$ viral
assemblies, specifically those of the viruses MVM (Minute Virus of
Mice), MSV (Maize Streak Virus) and AAV4 (Human Adeno Associate
Virus).

The GT model was developed to answer questions that concern
\emph{only} the influence of two quantities on the probability of each
type of
assembly tree.
We call these quantities the {\it geometric stability factor} and the
{\it
symmetry factor}. The higher these quantities, the higher the
probability.

The geometric stability factor is correlated with biochemical
stability and influenced by
assembly and disassembly energy thresholds and is
defined using the effect of geometric
constraints within monomers or between monomers. These constraints are
distances, angles and forces between the monomer residues.
These can be  obtained
either from X-ray or from cryo-electro-microscopic information on the
complete
viral shell.
More specifically, the final viral structure can be
viewed formally as the solution to a system of geometric constraints
that can be expressed as algebraic equations and inequalities.
For each internal vertex of an assembly tree, namely a subassembly,
the
geometric stability factor can be computed using quantifiable
properties
- such as extent of rigidity or algebraic complexity of
the configuration space of the subassembly.

\newpage

This factor is computed by analyzing the corresponding
subsystems of the given viral geometric constraint system.  It was
argued in \cite{bib:AgSi06} that the rigidity aspect of the
geometric stability factor can
be generically expressed purely using graph theory.
Some assembly trees can
never occur (have probability zero) since the subassemblies occurring
in them are unstable (their geometric stability factor is zero).  Such
assembly trees are {\em geometrically invalid.}

The symmetry factor is defined as follows.
The icosahedral group acts naturally 
 on the set of assembly trees for a particular viral polyhedron $P$, whose
 facets (representing viral monomers) are the leaves of the tree as 
illustrated in Figure \ref{pathways1}.
  Each orbit under this
 action is called an
{\it assembly pathway} and corresponds intuitively to a distinct type of 
assembly process for the viral capsid.  The symmetry factor
 is the number of assembly trees in the pathway divided by the total number
  of trees. We assume that each assembly tree is equally likely to occur.

In \cite{bib:icba}, the authors observed the following attractive
feature of the GT model of assembly.  The two separate factors -
geometric stability and symmetry - that influence the probability of the
occurrence of a particular assembly pathway can be analyzed largely
independently as follows.  An obvious, but crucial, observation
made in \cite{bib:icba} is that both the geometric stability factor
and geometric validity are invariants of the assembly pathway. That is,
they remain the same for any assembly tree in the same orbit under the
action of the icosahedral group.
Thus the probability of the occurrence of a
 pathway is roughly proportional to
some combination of the symmetric factor and the geometric
stability factor. Additionally, 
the ratio of the orbit sizes of two trees $\tau_1$ and $\tau_2$
could serve as a rough estimate of the the ratio of the probabilities
of the corresponding assembly pathways - {\em provided} that the
former ratio is  not cancelled out or reversed by the ratio of
the geometric stability factor of $\tau_1$ and $\tau_2$.  The paper
\cite{bonasith} formally proved that this kind of cancelling out
would not generally take place, at least for valid pathways, for the following
reasons.  First, it is shown in \cite{bonasith} that the
symmetry factor of a pathway increases with the depth of its
representative tree $\tau$.  More precisely, it was proved formally that the
size of the orbit of $\tau$ is bounded below by the depth of $\tau$.
Moreover, it is known from \cite{bib:AgSi06} that, {\em provided an
assembly tree is valid}, the geometric stability factor is non-zero
and generally increases with the depth of the tree (and this correlates with
biophysical observations).  Therefore, if the depth of $\tau_1$ is
greater than the depth of $\tau_2$, then {\em both} the symmetric
factor and the geometric stability factor of $\tau_1$ will generally be larger
than the corresponding factors of $\tau_2$.  


\subsection{Contributions and Related Work}

Based on the observations in the last section, the paper
\cite{bonasith} posed problems intended to isolate and clarify the
influence of the symmetry factor on the probability of the occurrence
of a given assembly pathway.  Two specific problems were the
following.
\begin{enumerate}
\item [(i)] Enumerate the valid assembly pathways of an icosahedrally
symmetric polyhedron. More precisely, the problem is to determine the
number of such assembly pathways of each orbit size. 

\item [(ii)] Characterize and algorithmically recognize the set
of assembly trees fixed by a given subgroup of the icosahedral
group. The characterization problem is a step toward the solution of
the enumeration problem (i).  Algorithmic recognition of the
group elements that fix a given assembly tree (the stabilizer of 
the tree in the given group) 
directly determines whether the given assembly
tree has a given orbit size.
\end{enumerate}

{\bf Remark.}
{\em In this paper, we answer the above questions for general assembly trees,
that is, {\it we drop the condition of validity}. Furthermore, in this
paper, the geometric stability factor will be ignored,  and thus 
``probability'' will refer to tbe symmetry factor only.  As mentioned earlier,
\cite{bib:AgSi06,bib:icba} 
show that validity of assembly trees of a polyhedron
$P$ is not only invariant under the action of its symmetry group, but 
can also be captured by simple graph-theoretic properties such as
generalized notions of connectivity for the graph  constructed from
 the vertices and edges of $P$.  We expect that the techniques developed in
this paper will help in answering the above questions in the presence of
the validity condition as well.}

\vskip 0.3 cm 
For Problem (i), we develop an enumeration method using generating
functions and M\" obius inversion. For the algorithm in Problem (ii), we
provide a simple permutation group algorithm and an associated data
structure.  The results of this paper work not just for the
icosahedral group, but also for any finite group $G$ acting simply on a set
$X$. Indeed, if $G$ is a finite group acting on a set $X$, then there
is a natural induced action of $G$ on the set ${\mathcal T}_X$ of {\it
assembly trees}. These are formally defined as rooted trees $\tau$ whose
non-leaf vertices have at least two children and whose leaves are
bijectively 
labeled by $X$. If $G$ acts simply on $X$, then $|X| := |X_n| = n \,
|G|$, where $n$ is the number of $G$-orbits in $X$.

Concerning Problem (i), P\'olya theory gives a convenient method
for counting orbits under a permutation group action. 
However, because
of the complexity of the cycle index in our situation, we were not
able to apply P\'olya theory to Problem (i).
Similarly, the methods used in \cite{bib:klin}
 for enumerating labeled graphs under a
group action (as opposed to rooted labeled trees), did not seem to
apply. Our generating function method, on the other hand, finds an
explicit formula (Theorem~\ref{orbit}) for the number of orbits of
each possible size in the action of $G$ on the set ${\mathcal T}_{X_n}$ of
assembly trees, for every $n$.  This leads to a formula for the
probability of occurrence of a given assembly pathway
(Corollary~\ref{prob}).  To apply these formulas
it is necessary to know the number of assembly trees fixed by
each given subgroup of G. A generating function formula for this number of
fixed assembly trees is given in Theorem 16 of Section 5. For the proof
of Theorem 16 is is necessary to characterize the set of such fixed
assembly trees. This is done in Theorem 9 of Section 4.

Concerning Problem (ii), algorithms for permutation groups have been
well-studied (see for example \cite{bib:permalgo}), and 
algorithms for tree isomorphism and automorphism are well known 
\cite{bib:treecompare,bib:treealgobook}.  
Moreover, the structure of the automorphism groups of rooted, labeled 
trees have been studied \cite{bib:automorph,bib:polyaidentity}. 
%
%
However,
 we have not encountered an algorithm in the literature for deciding whether a
given permutation group element fixes a given rooted, labeled tree;
and thereby finds the stabilizer of that tree in the given group $G$. 
In Section 3.2 of this paper, we
provide a simple and intuitive algorithm that is easy to
implement, runs in
linear time and operates in place on the input, without the use of  extra
scratch memory.

If one is only interested in approximate and asymptotic estimates for
Problem (i), such as in viruses with large T-numbers, a possible avenue is 
to use the results of
\cite{bib:compton, bib:woods} that estimate the asymptotic
probabilities of logic properties on finite structures, especially
trees. There are significant roadblocks, however, to applying these
results to our problem.  These are mentioned in the open problem
section at the end of this paper. Finally, there is a rich literature
on the enumeration of construction sequences of symmetric polyhedra
and their underlying graphs \cite{Brinkmann,Dutour,Fowler}. Whereas
these studies focus on enumerating construction sequences of {\em
different} polyhedra with a given number of facets, our goal - of
counting and characterizing assembly tree orbits - is geared towards
enumerating construction sequences of a {\em single} polyhedron for
any given number of facets.

\section{Preliminaries on Assembly Pathways}
\label{sec:definition}

All groups, graphs, and label sets in this paper are assumed to be finite.
A {\it rooted tree} is a tree with a designated vertex, called the
{\it root}.  We will use standard terminology such as {\it adjacent},
{\it child}, {\it parent}, {\it descendent}, {\it ancestor}, {\it
leaf}, {\it subtree rooted at}, {\it root of the subtree}, and so on.  For
our purposes, a rooted tree is called a {\it labeled} tree if the
leaves are bijectively labeled by the elements of a set $X$, and an
{\it internal} (non-leaf) vertex $v$ is labeled by the set of
leaf-labels of the subtree rooted at $v$.  {\em We identify each vertex in
a tree with its label.}  Let $\tau$ and $\tau'$ be two rooted trees
labeled in  the same set $X$.
Then $\tau$ and $\tau'$ 
are said to be {\it isomorphic} if there is a bijection - the {\it
isomorphism} - $f$ between the vertices of $\tau$ and $\tau'$ that
preserves adjacency and the root. That is, all the following hold.
\begin{itemize}
\item  $(u,v)$ is an edge in
$\tau$ if and only if $(f(u),f(v))$ is an edge in $\tau'$, 
\item 
$f(r) = r'$, where $r$ and $r'$ are the roots of $\tau$ and $\tau'$
respectively  In this case, we say $\tau \approx \tau'$ and also
$f(\tau) = \tau'$. 
\end{itemize}
  An {\it automorphism} $f$ of
$\tau$ is an isomorphism of $\tau$ into itself: it ensures that
$(u,v)$ is an edge in $\tau$ if and only if $(f(u),f(v))$ is also an
edge in $\tau$.  In this case $f(\tau) = \tau$. \m

\begin{figure}[!ht]
\vskip 2mm
\begin{center}
\includegraphics[width=4.5in]{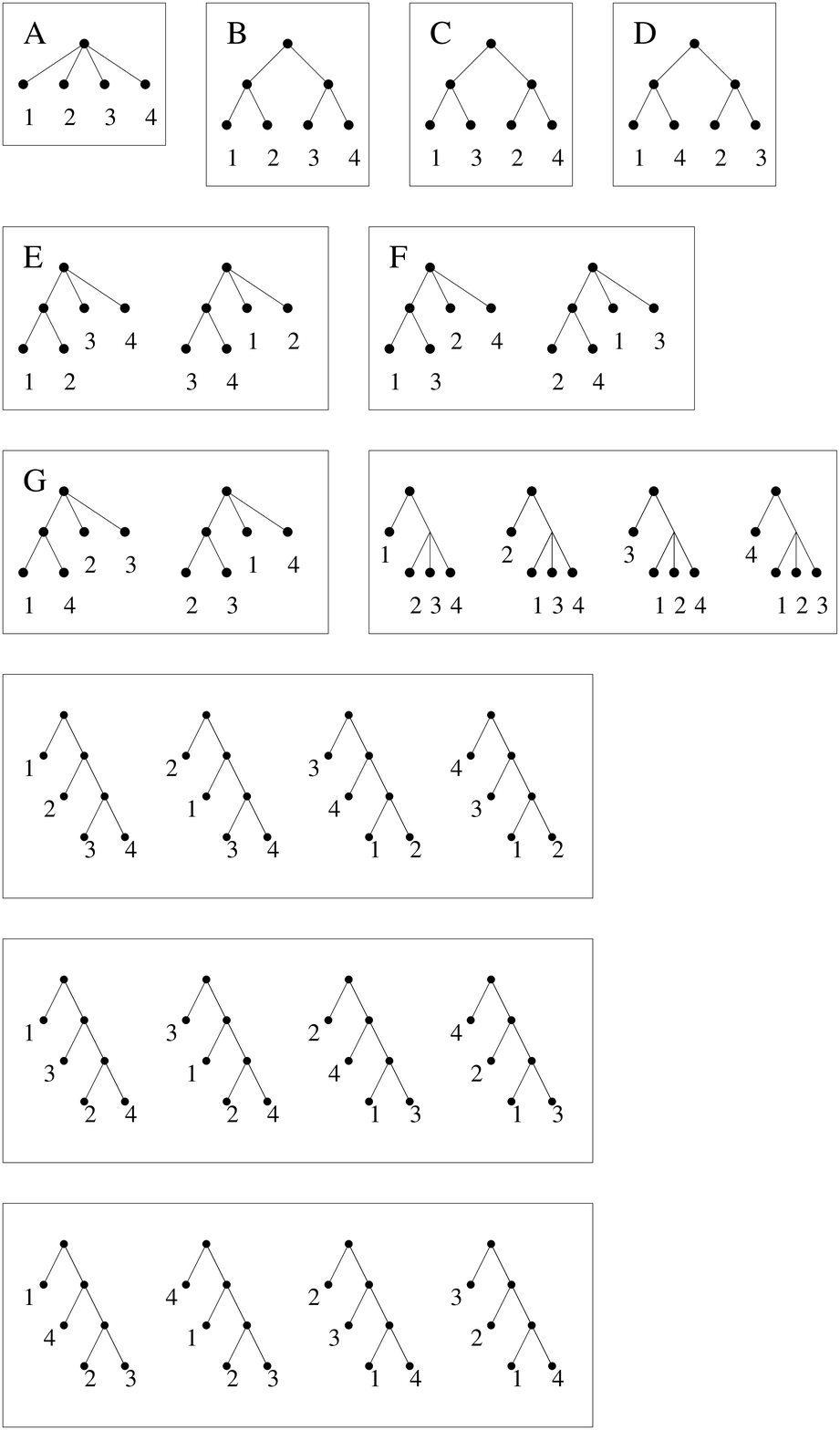}
\end{center}
\caption{Klein 4-group acting on $\mathcal{T}_4$.} 
\label {trees4}
\end{figure}

A rooted tree for which each internal vertex has at least two
children and whose leaves are labeled with elements of $X$ is called
an {\it assembly tree for} $X$.  The $26$ assembly trees with four
leaves, labeled in the set $X = \{1,2,3,4\}$ are shown in
Figure~\ref{trees4}.  \m

Let $G$ be a group acting on a set $X$.  
The action of $G$ on $X$ induces a natural action of $G$ on 
the power set of $X$ 
and thereby on the set of vertices (vertex labels) of 
$\mathcal{T}_X$ of assembly trees for $X$.  
If $g \in G$ and $\tau\in \mathcal{T}_X$,
then define the tree $g(\tau)$ as the unique assembly tree 
whose set of vertex labels (including the labels of internal vertices) 
is $\{g(v): v\in \tau\}$. 
This tree $g(\tau)$ is clearly isomorphic to $\tau$ via $g$.
This induces an action of $G$ on  $\mathcal{T}_X$.
Each orbit of this action of $G$ on 
$\mathcal{T}_X$ consists of isomorphic 
trees and is
called an {\it assembly pathway for} $(G,X)$.

\begin{example} Klein 4-group acting on $\mathcal{T}_4$. 
\label{klein}\end{example}

Consider the Klein 4-group $G = \Z _2 \oplus \Z _2$ acting
on the set $X = \{1,2,3,4\}$. Writing $G$ as a group of 
permutations in cycle
notation, this action is
$$G = \{(1)(2)(3)(4), \, (1\;2)(3\;4), \, (1\;3)(2\;4), \,
(1\;4)(2\;3)\}.$$  For
this example there are exactly 11 assembly pathways, which are
indicated in Figure~\ref{trees4} by boxes around the orbits.  There
are four assembly pathways of size one, i.e., with one assembly tree in the
orbit, three assembly pathways of size two, and four assembly pathways of size
four. \B


An assembly tree $\tau$ is said to be
{\it fixed} by an element $g\in G$ if $g$ is an automorphism of $\tau$, 
that is, $g(\tau) = \tau$. See Figure \ref{pathways1} for an
illustration. For any
subgroup $H$ of $G$, let $\overline{t}(H) := \overline{t_X}(H)$ denote
the number of trees in $\mathcal{T}_X$ that are fixed by all elements of
$H$ and by no other elements of $G$. In other words,
\begin{equation}
\overline{t}_X(H) =
 |\{ \tau \in \mathcal{T}_X \, | \, stab_G(\tau) = H \}|. 
\label{orbitsize}
\end{equation}
Here  $stab_G(\tau) :=
\{ g \in G \; | \; g(\tau)=\tau\}$ is called the {\em stabilizer} of
$\tau$ in $G$. In other words,  $stab_G(\tau)$ is the set of all elements 
in $G$ that fix $\tau$. It is easy to prove that  $stab_G(\tau)$ is
a subgroup of $G$.  
\m

In fact, we will see in Section \ref{enumerating} that it is more natural
to find $t_X(H)$, i.e., the number of trees in $\mathcal{T}_X$
that are fixed by  a subgroup $H$ of $G$. These may include trees
that are fixed by larger subgroups $H'$ such that $H \le H' \le G$.
As the following theorem shows, the desired quantities 
$\overline{t}_X(H)$ can then be computed from the numbers $t_{X}(H)$
 using M\"obius inversion on
the lattice of subgroups of $G$.

\begin{theorem} \label{mobius} Let $G$ be a group acting on a set $X$.
If $H$ is a subgroup of $G$, then
$$\overline{t_X}(H) = \sum_{H\leq K \leq G} \; \mu (H,K)\; t_X(K),$$
where $\mu$ is the M\"obius function for the lattice of subgroups of $G$.
\end{theorem}

\proof Clearly $t_X(H) = \sum_{H\leq K \leq G} \overline{t_X}(K)$.
The theorem follows from the standard M\"obius inversion formula
\cite{VW} (page 333).
\qed \B

The index of a subgroup $H$ in $G$ is the number of left
(equivalently, right), cosets of
$H$ in $G$, and is 
denoted by $(G:H)$.  By Lagrange's Theorem, this index equals $|G|/|H|$.

\begin{theorem} \label{orbit}
The number of trees in any assembly pathway for $(G,X)$ divides $|G|$.
If $m$ divides $|G|$, then the number $N(m)$ of assembly pathways of
size $m$ is
 $$N(m) = \frac{1}{m} \; \sum_{H\leq G \,: \, (G:H) = m} \;
\overline{t}(H).$$
\end{theorem}

\proof  It is a standard consequence of Lagrange's Theorem that, for
any assembly tree $\tau$, the equality  
$$|G| = |O(\tau)| \cdot |stab(\tau)|$$ holds,   where $O(\tau)$ is the 
orbit of $\tau$.
This immediately implies the first statement of the theorem.

Let
$$\delta(\tau) = \begin{cases} 1 &\text{if $(G: stab(\tau)) = m$} \\
0 & \text{otherwise} \end{cases} \quad
= \begin{cases} 1 &\text{if $|O(\tau)|= m$} \\
0 & \text{otherwise}.\end{cases}$$
Now count, in two ways, the number of pairs $(H,\tau)$ where $\tau$ is an
assembly tree, $H\leq G$ is the stabilizer of $\tau$, and $(G:H) = m$:
$$\sum_{H\leq G \,: \, (G:H) = m} \; \overline{t}(H) = 
\sum_{\tau \in {\mathcal T}_X} \delta (\tau) = m \, N(m).$$
Indeed, to justify the first equality, note that for a fixed subgroup
$H$ that has index $m$ in $G$, exactly $\bar{t}(H)$ trees $\tau$ will
satisfy $\delta(\tau)=1$. To justify the second equality, note that $\delta(\tau)=1$
if and only if $\tau$ is one of $m$ elements of an $m$-element pathway, and
there are $N(m)$ such pathways.
\qed \m

\begin{example} Klein 4-group acting on ${\mathcal T}_4$ (continued). 
\end{example} 
Theorem~\ref{orbit}, applied to our previous example of $\Z _2 \oplus \Z _2$
acting simply on $\{1,2,3,4\}$, states that the size of an assembly
pathway must be $1$, $2$ or $4$, since it must be 
a divisor of $4 = |\Z _2 \oplus \Z _2|$.
To find number of pathways of each size, note that
$G$ has three subgroups of order 2, namely
$$\begin{aligned}
K_1 &= \{ \,(1)(2)(3)(4), (1\;2)(3\;4) \, \}, \\
K_2 &= \{ \,(1)(2)(3)(4), (1\;3)(2\;4) \, \}, \\
K_3 &= \{ \,(1)(2)(3)(4), (1\;4)(2\;3) \, \},
\end{aligned}$$
and that
$$\begin{aligned}
& \overline{t}(G) = 4, \\
& \overline{t}(K_1) = \overline{t}(K_2) =\overline{t}(K_3) = 2, \\
& \overline{t}(K_0) = 16,
\end{aligned}
$$
where $K_0$ denotes the trivial subgroup of order 1.
The assembly trees in ${\mathcal T}_X$ that are fixed by all elements
of $G$ are shown in Figure~\ref{trees4}, $A, B, C, D$.  For $i =
1,2,3$, those assembly trees in ${\mathcal T}_X$ that are fixed by all
elements of $K_i$ and by no other elements of $G$ are are shown in
Figure~\ref{trees4}, $E, F, G$, respectively. The remaining 16 assembly
trees in Figure~\ref{trees4} are fixed by no elements of $G$ except
the identity.  Therefore, according to Theorem~\ref{orbit},
the number of pathways of size 1, 2 and 4 are, respectively,
$$\begin{aligned}
\overline{t}(G) &= 4, \\
\frac12 \, \left( \, \overline{t}(K_1) + \overline{t}(K_2) +
\overline{t}(K_3)\, \right) = \frac12 \, (2+2+2) &= 3, \\
\frac14 \, \overline{t}(K_0) &= 4.
\end{aligned}
$$
A general formula for $\overline{t}(H)$ is the subject of
Sections~\ref{characterizing} and \ref{enumerating}. \B

The set ${\mathcal T}_X$ of assembly trees can be made the sample space of
a probability space $({\mathcal T}_X, p)$ by assuming that each assembly
tree $\tau \in {\mathcal T}_X$ is equally likely, i.e., $p(\tau) = 1/|{\mathcal
T}_X|$. Clearly, if $O$ is an assembly pathway, then 
$$p(O) = \frac{|O|}{|{\mathcal T}_X|}.$$  The following result
follows immediately from Theorem~\ref{orbit}.

\begin{cor} \label{prob} If $G$ acts on the set $X$ and $m$ divides
$|G|$, then, with notation as in Theorem~\ref{orbit}, there are exactly
$N(m)$ assembly pathways with probability $\frac{m}{|{\mathcal T}_X|}$,
and no other values can occur as the probability of an assembly
pathway.
\end{cor}

\begin{example}  Klein 4-group acting on ${\mathcal T}_4$ (continued). 
\end{example} 
Again, for our example of $\Z _2 \oplus \Z _2$ acting simply on
$\{1,2,3,4\}$, application of Corollary~\ref{prob} gives
$$\begin{aligned}  4& \quad
\hbox{pathways with probability}  \quad \frac{1}{26}, \\
3& \quad
\hbox{pathways with probability}  \quad\frac{1}{13}, \\ 
 4& \quad
\hbox{pathways with probability} \quad
\frac{2}{13}.
\end{aligned}$$

\section{Algorithm for determining the stabilizer of an assembly 
tree in a given finite group}
\label{sec:algorithm}

The algorithm in this section takes as input a finite permutation
group $G$ acting on a finite set $X$ and an assembly tree $\tau \in
\mathcal{T}_X$, and finds the stabilizer $stab_G(\tau)$.  The idea
behind the algorithm is encapsulated by the following proposition,
whose proof follows directly from definitions given in Section 2.  As defined
in Section~\ref{sec:definition}, the action of the permutation group
 $G$ on $X$
induces a natural action of $G$ on $\mathcal{T}_X$.
\newpage

\begin{proposition}
\label{algo-correct}
Let the finite permutation group $G$  act on  a finite set $X$. 

\begin{enumerate}
\item
Let $R$ be any set of elements of $G$ that fix $\tau$, and let 
$C$ be any set of elements of $G$ that do not fix $\tau$.
Then $\langle R\rangle$, the group generated by the elements of $R$,
is a subgroup of $stab_G(\tau)$ and 
$\bigcup\limits_{c\in C} c\langle R\rangle$, the union of the 
left cosets of $\langle R \rangle$
given by $C$, has an empty intersection with  
$stab_G(\tau)$.
\item
An element $g \in G$ fixes $\tau$ if and only if
for every vertex $v\in \tau$ with children $c_1(v),\ldots, c_k(v)$,
the vertices $g(c_1(v)),\ldots ,g(c_k(v))$ have a common parent in $\tau$,
and this parent is $g(v)$.
\end{enumerate}

\end{proposition}

\subsection{Input and Data Structures}
In this subsection, we give the detailed setup 
for our stabilizer-finding algorithm. 
The input of this algorithm is the set of elements $g$
of the finite permutation group $G$ acting on the finite set $X$
and an assembly tree $\tau \in \mathcal{T}_X$. 

We use a  tree data structure, where 
each vertex has a child pointer to each of its children and 
a parent pointer to its parent.
The root (and the tree $\tau$ itself) can be accessed by the root pointer.
Furthermore, each permutation $g\in G$ on $X$ is input as a set of 
$g$-pointers on the leaves.
That is,  a leaf labeled $u$ has a $g$-pointer to the leaf labeled $g(u)$. 
However, the labels are not explicitly stored except at the leaves.
This is a common data structure used in permutation group algorithms
\cite{bib:permalgo}. 


\subsection{Algorithms}

The first algorithm computes $stab_G(\tau)$, and it  uses
the second and third algorithms for determining whether a permutation $g$
fixes $\tau$.
The correctness of the first  algorithm follows  directly 
from Proposition \ref{algo-correct} (1), assuming the correctness of 
the latter algorithms.
The last algorithm is recursive and operates in place with no extra scratch
space. 
For each vertex $v$, working bottom up, it 
efficiently checks  whether the image 
$g(v)$ is in $\tau$. The correctness follows directly from 
Proposition \ref{algo-correct} (2). 

\medskip\noindent
{\bf Algorithm Stabilizer}\\
{\it Input:} assembly-tree $\tau \in \mathcal{T}_X$; permutation group,
$G$ \\ 
{\it Output:} generating set $R_\tau$ s.t. the group $\langle R_\tau \rangle$ 
generated by $R_\tau$ is exactly $stab_G(\tau)$.\\

\noindent
$R := \{id\}$ ({\tt currently known partial generating set of} 
$stab_G{(\tau)}$) \\
$C_R := \emptyset$ ({\tt distinct left coset representatives of} 
     $\langle R \rangle$ \\
$\ ~~~~~~~~~~${\tt that are currently known to not fix} $\tau$)\\
$U :=  G$ ({\tt currently undecided elements of} $G$)\\
$\ ~~${\tt do until} $U = \emptyset$\\
$\ ~~${\tt let}  $g\in U$ \\
$\ ~~~~$ {\tt if} {\bf Fixes}$(g,\tau)$\\ 
$\ ~~~~~~${\tt then} $R := R\cup  \{g\}$; {\tt retain in} 
$C_R$ {\tt at most one\\ 
$\ ~~~~~~~~~~~~~~$representative from any left coset of} $\langle R\rangle$\\
$\ ~~~~~~${\tt else} $C_R := C_R \cup \{g\};$\\
$\ ~~~~$$U:=(U \setminus(\langle R \rangle \bigcup\limits_{c\in C_R}
 c\langle R\rangle)$ \\
$\ ~~~~${\tt fi}\\
$\ ~~${\tt od}\\ 
{\tt return} $R_\tau := R.$\\

\bigskip\noindent
{\bf Algorithm Fixes}\\
{\it Input:} assembly-tree $\tau \in \mathcal{T}_X$; permutation 
$g$ acting on $X$,\\ 
{\it Output:} 
``true'' if $g$ fixes $\tau;$ ``false'' otherwise.\\

\noindent
{\tt if} {\bf LocateImage}$(g,\tau,root(\tau)) = root(\tau)$\\ 
{\tt then return true}\\
{\tt else return false.}

\bigskip\noindent
{\bf Algorithm LocateImage}\\
{\it Input:} assembly-tree $\tau \in \mathcal{T}_X$; permutation 
$g$ acting on $X$; child pointer to a vertex $v \in \tau$ (root pointer
if $v$ is the root of $\tau$).\\ 
{\it Output:} 
a parent pointer to the vertex $g(v) \in \tau$ - if it exists -
such that $g$ is the isomorphism mapping the subtree of $\tau$ 
rooted at $v$ to the subtree of $\tau$ rooted at $g(v)$;
if such a vertex $g(v)$ does not exist in $\tau$, returns null.\\

\noindent
{\tt if} $v$ {\tt is a leaf of} $\tau$ {\tt (null child pointer)} \\
{\tt then return}
$g(v)$ {\tt (follow} $g${\tt -pointer)}\\
{\tt else}\\
$\ ~~$
{\tt let} $c_1,\ldots,c_k$ {\tt be the children of} $v$;\\
$\ ~~$
{\tt if} {\tt parent}({\bf LocateImage}$(g,\tau,c_1)$) = \\
$\ ~~~~$ 
{\tt parent}({\bf LocateImage}$(g,\tau,c_2)$) = \ldots\\
$\ ~~~~$
{\tt parent}({\bf LocateImage}$(g,\tau,c_k)$) =: $w$ \\
$\ ~~$
{\tt then return} $w$\\ 
$\ ~~$
{\tt else return null.}

\subsection{Complexity}

{\bf Algorithm LocateImage} follows each pointer (child, $g$, parent)
exactly once as is illustrated by the example shown in Figures 6-9,
 and does only constant time
operations between pointer accesses. Hence it takes at most $O(|X|)$
time. It operates in place and does not require any extra scratch
space.  {\bf Algorithm Stabilizer}, 
in the worst case, can be a brute
force algorithm that simply runs through all the elements of $G$
instead of maintaining efficient representations of $R$, $C_R$ and $U$.  
In this case, it takes no more than $O(|G||X|)$ time.

However, readers familiar with Sim's method for representing 
permutation groups \cite{bib:permalgo} using so-called strong generating sets 
and  Cayley graphs may appreciate the
following remarks. Instead of specifying the input of our problem as
we have done, we may assume that a Cayley graph is
input, which uses a strong generating set of $G$. 
With this  input representation, the time complexity of our
algorithms can be significantly 
further optimized, the level of optimization depending on 
properties of the group $G$.

\subsection{Example}

\begin{figure}
\psfrag{1}{$1$} \psfrag{2}{$2$} \psfrag{3}{$3$} \psfrag{4}{$4$}
\psfrag{c12}{$\{1,2\}$} \psfrag{c1234}{$\{1,2,3,4\}$}
\psfrag{g1234}{$g = (1,2)(3,4)$} \centerline{\epsfig{file =
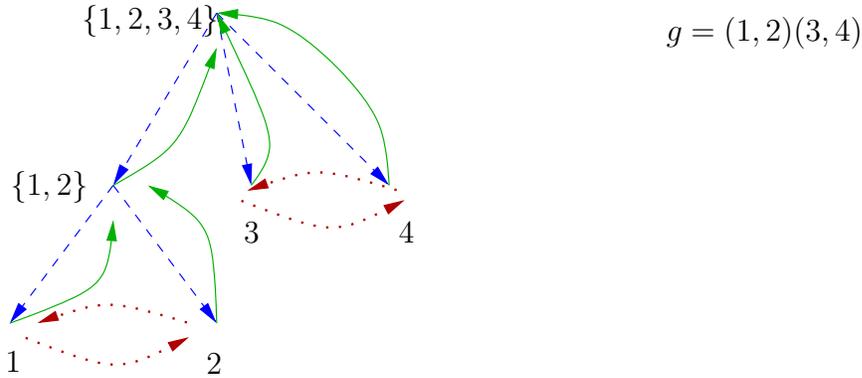, height = 2in}}
\caption{{\bf LocateImage} is called on the root of the assembly tree
$\tau$ shown on the left, for the permutation $g$ shown on right. This
permutation $g$ fixes $\tau$.  The data structure representing the tree
consists of the child (blue, dashed) and parent (green, solid)
pointers and the data structure representing $g$ - via its action on
the leaf-label set $X$ - consists of $g$-pointers (red, dotted).}
\label{fig:automorph}
\end{figure}

\begin{figure}
\begin{center}
\psfrag{1}{$1$}
\psfrag{2}{$2$}
\psfrag{3}{$3$}
\psfrag{4}{$4$}
\psfrag{c12}{$\{1,2\}$}
      \psfrag{c1234}{$\{1,2,3,4\}$}
      \psfrag{l1}{$\ $}
      \psfrag{l2}{$\ $}
      \psfrag{l3}{$g(1)=2$)}
      \psfrag{l4}{LocateImage($1$)}
      \psfrag{l5}{LocateImage($\{1,2\}$)}
      \psfrag{l6}{LocateImage(root)}
\epsfig{file = 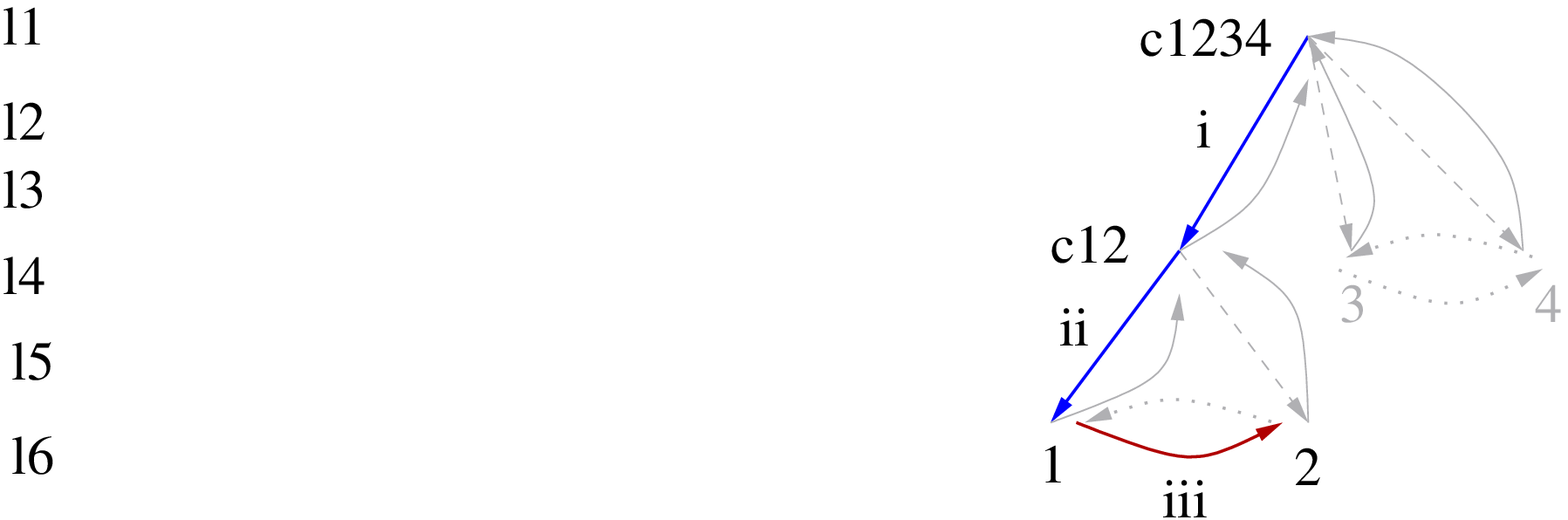, height = 1.5in}\\
      \psfrag{l1}{$\ $}
      \psfrag{l2}{$\ $}
      \psfrag{l3}{$g(2)=1$}
      \psfrag{l4}{$g(1)=2$; parent($g(1)=\{1,2\}$)}
      \psfrag{l5}{LocateImage($\{1,2\}$)}
      \psfrag{l6}{LocateImage(root)}
\epsfig{file = 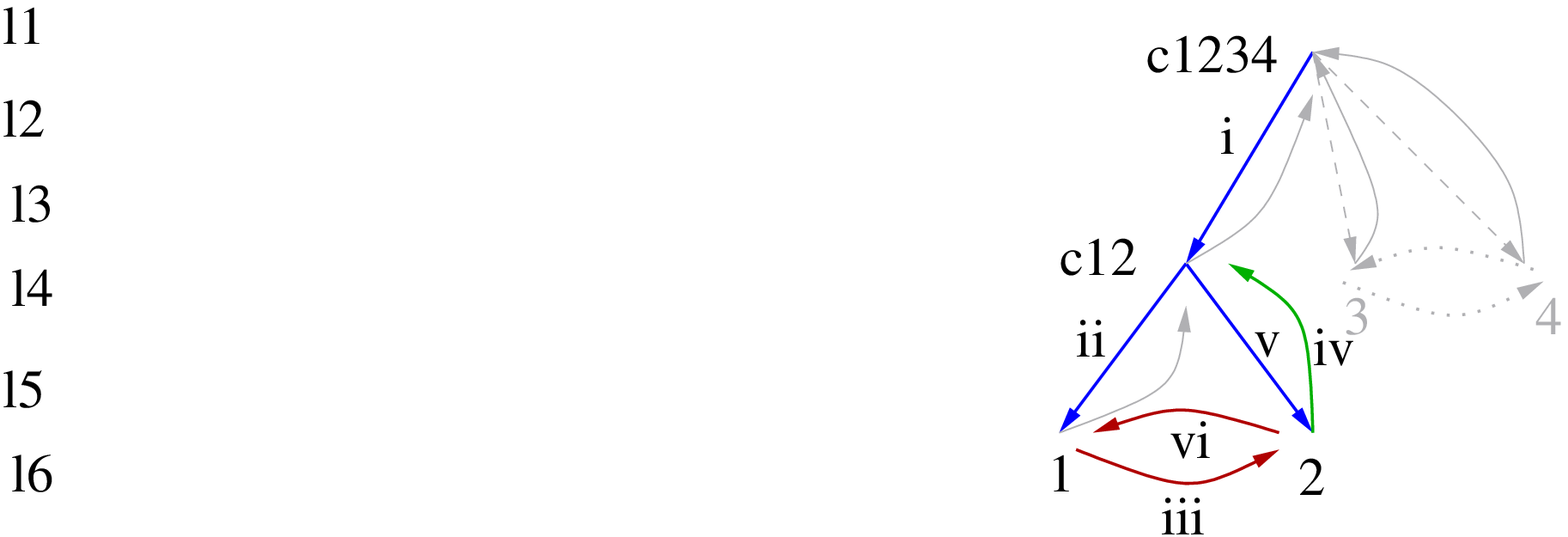, height = 1.5in}\\
      \psfrag{l1}{$\ $}
      \psfrag{l2}{LocateImage($3$)}
      \psfrag{l3}{$g(\{1,2\})=\{1,2\}$}
      \psfrag{l4}{$g(2)=1$; parent($g(2)=\{1,2\}$)}
      \psfrag{l5}{$g(1)=2$; parent($g(1)=\{1,2\}$)}
      \psfrag{l6}{LocateImage(root)}
\epsfig{file = 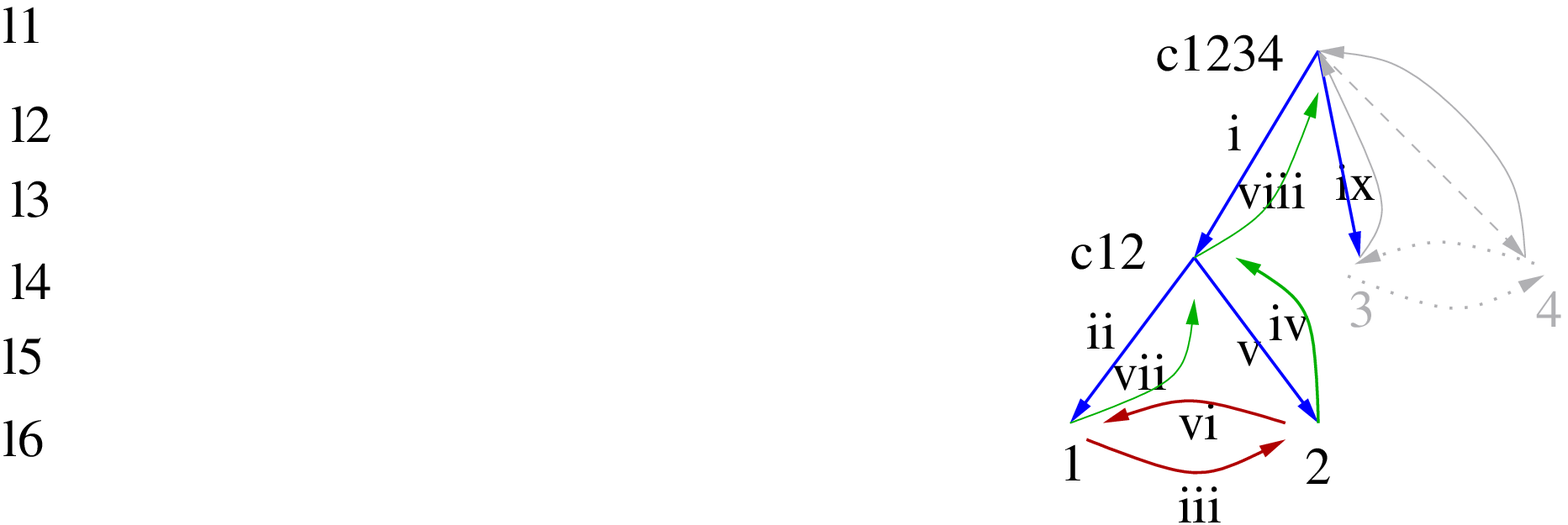, height = 1.5in}\\
      \psfrag{l1}{$g$(root)=root; Fixes returns True}
      \psfrag{l2}{$g(4)=3$; parent($g(4)=$ root}
      \psfrag{l3}{$g(3)=4$; parent($g(3)=$ root}
      \psfrag{l4}{$g(\{1,2\})=\{1,2\}$; parent($g(\{1,2\})$)=root}
      \psfrag{l5}{$g(2)=1$; parent($g(2)=\{1,2\}$)}
      \psfrag{l6}{$g(1)=2$; parent($g(1)=\{1,2\}$)}
$\ $\\
\epsfig{file = 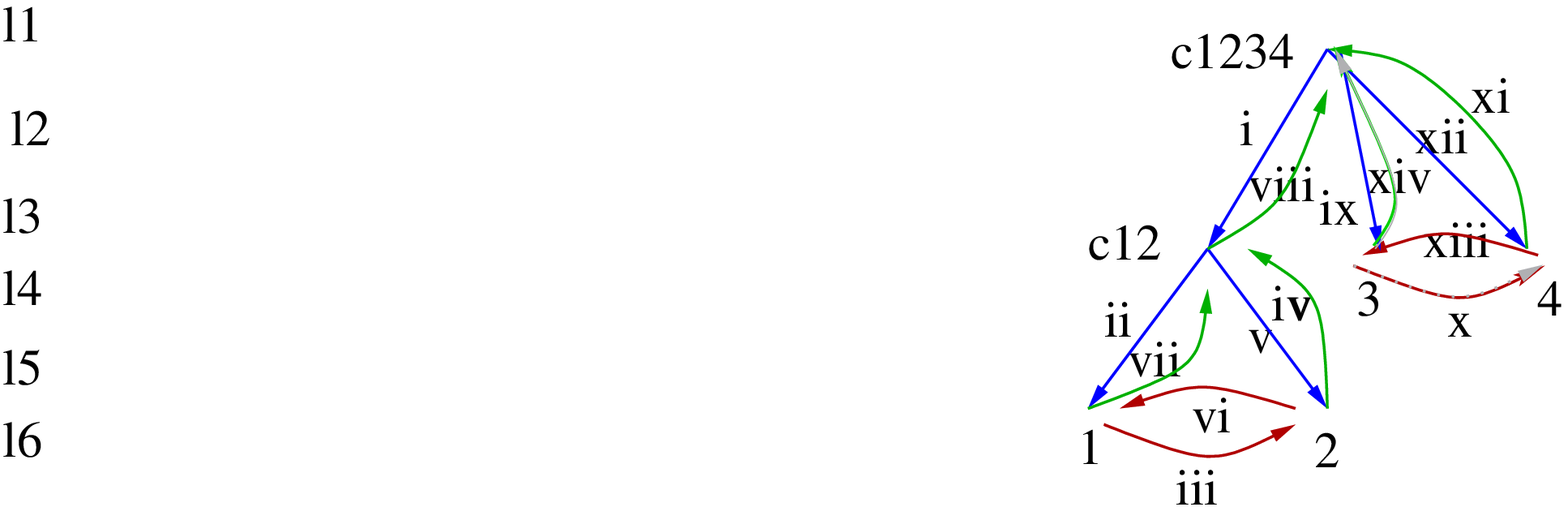, height = 1.5in}
\end{center}
\caption{A successful run of {\bf LocateImage},
 where, for all vertices
$v$ in $\tau$, the image $g(v)$ is established to be in the
assembly-tree $\tau$ shown in Figure \ref{fig:automorph}.  On the
right are pointers traversed so far, in traversed order.  On the left
is the current recursion stack of {\bf LocateImage} calls (first call
at the bottom), together with those vertices $v$ ($\in X$ or
$\subseteq X$) for which $g(v)$ has been established to be in $\tau$,
showing that $g$ is an isomorphism between the two subtrees of $\tau$
rooted at $v$ and at $g(v)$, respectively.} 
\label{fig:success}
\end{figure}

\begin{figure}
\psfrag{1}{$1$}
\psfrag{2}{$2$}
\psfrag{3}{$3$}
\psfrag{4}{$4$}
\psfrag{c12}{$\{1,2\}$}
      \psfrag{c1234}{$\{1,2,3,4\}$}
      \psfrag{g1423}{$g = (1,4)(2,3)$}
\centerline{\epsfig{file = 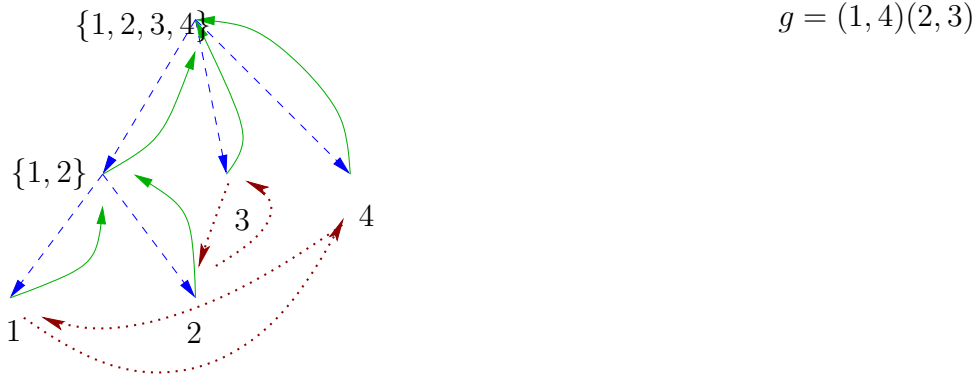, height = 2in}}
\caption{{\bf LocateImage} is called on the root of the  
same tree $\tau$ as in Figure \ref{fig:automorph}, 
but for a different permutation $g$ shown on right. In this case 
 $g$ does {\sl not} fix $\tau$. }
\label{fig:nonautomorph}
\end{figure}

\begin{figure}
\psfrag{1}{$1$}
\psfrag{2}{$2$}
\psfrag{3}{$3$}
\psfrag{4}{$4$}
\psfrag{c12}{$\{1,2\}$}
      \psfrag{c1234}{$\{1,2,3,4\}$}
      \psfrag{L1}{g(root)=null; Fixes returns False}
      \psfrag{L2}{parent($g(3)$)$\not=$ parent($g\{1,2\}$)}
      \psfrag{L3}{$g(3)=2$; parent($g(3)= \{1,2\}$}
      \psfrag{L4}{$g(\{1,2\})=\{1,2\}$; parent($g(\{1,2\})$)=root}
      \psfrag{L5}{$g(2)=3$; parent($g(2)$)=root)}
      \psfrag{L6}{$g(1)=4$; parent($g(1)$)=root)}
\centerline{
\epsfig{file = 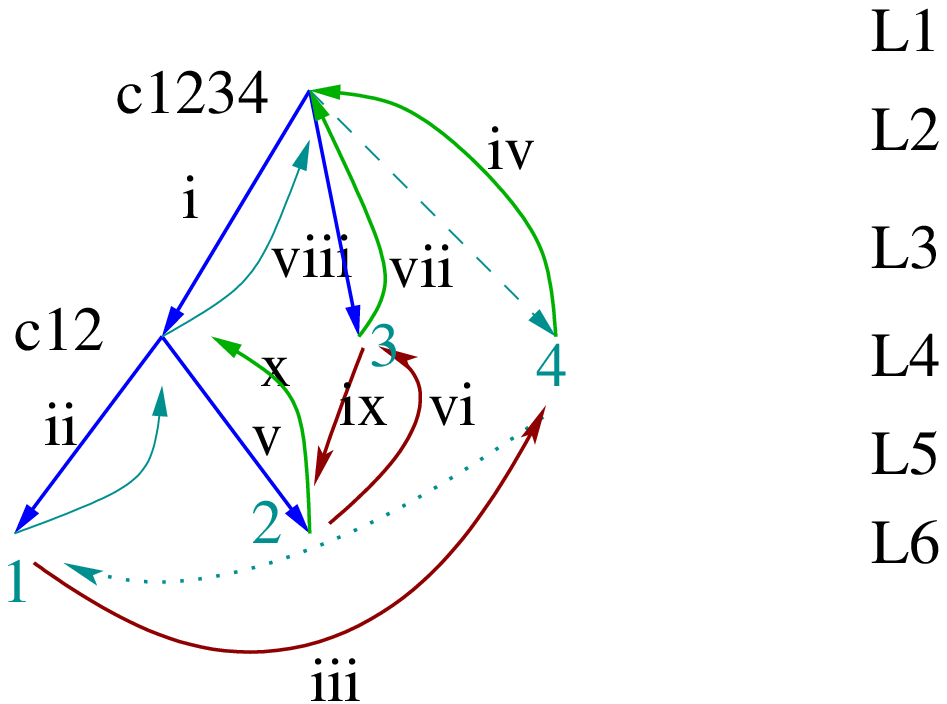, height = 2in}\\
}
\caption{An unsuccessful run of {\bf LocateImage}. 
Since $\{1,2\}$,
$3$ and $4$ are the children of the root, when {\bf LocateImage}(root)
is called, it checks if their images under $g$ have the same parent.
So recursive calls are to {\bf LocateImage}$(\{1,2\})$, to {\bf
LocateImage}$(3)$, and to {\bf LocateImage}$(4)$.  {\bf
LocateImage}$(\{1,2\})$ returns the root of $\tau$ as a candidate for
$g(\{1,2\})$.  But {\bf LocateImage}$(3)$ returns $\{1,2\}$ because
the parent of $g(3) = 2$ is $\{1,2\}$.  Hence
{\bf LocateImage}(root) returns null and {\bf Fixes} returns `false.'}
\label{fig:nonsuccess}
\end{figure}

The two examples shown in Figures 6-9
 illustrate the algorithms {\bf LocateImage} and
{\bf Fixes}. Figure \ref{fig:automorph} for the first example shows
the assembly tree $\tau$ and the associated data structures, as well
as the group element $g$. Figure \ref{fig:success} shows a run of {\bf
LocateImage} applied to $\tau, g$ at $root(\tau)$.  The algorithm
establishes that the given permutation $g$ fixes the given assembly
tree $\tau$, whereby {\bf Fixes} returns `true.'  The second example,
in Figure \ref{fig:nonautomorph}, uses the same assembly tree $\tau$,
but a different permutation $g$.  The run of {\bf LocateImage} in
Figure \ref{fig:nonsuccess} is unsuccessful, whereby {\bf Fixes}
returns `false.'

\section{Block systems and Fixed Assembly Trees} 
\label{characterizing}

The formulas in Section \ref{sec:definition} 
for the number of orbits of each size and
for the orbit sizes or pathway probabilities (Theorem~\ref{orbit} and
Corollary~\ref{prob}) depend on the number of assembly trees fixed by
a group.  A formula for the number of such fixed trees is
the subject of this and the next section.

Recall that an assembly tree $\tau$ is fixed by a group $G$ acting on $X$
if $g(\tau) = \tau$ for all $g\in G$.  Two main results of this
section (Corollary~\ref{a-tree} and Procedure~\ref{const}) provide a
recursive procedure for constructing all trees in ${\mathcal T}_X$ that
are fixed by $G$. This leads, in the next section, to a generating
function for the number of such fixed trees.  The results in this
section depend on a characterization (Theorem~\ref{CC}) of 
block systems arising from a group acting on a set.  \m

For a group $G$ acting on set $X$, a {\it block} is a subset
$B\subseteq X$ such that for each $g\in G$, either $g(B) = B$ or $g(B)
\cap B = \emptyset$. A {\it block system} is a partition of $X$ into
blocks.  A block system $\bf B$ will be said to be {\it
compatible} with the group action if $g(B)\in {\bf B}$ for all $g \in
G$ and $B \in {\bf B}$. A characterization of complete block systems 
(Theorem~\ref{CC}) is
relevant to the understanding of fixed assembly trees because of the
following result.  Let $\tau$ be any assembly tree in ${\mathcal T}_X$.  For
any vertex $v$ of $\tau$, recall that $v$ is identified with 
and labeled by its set of descendent leaf-labels.
Thus the  set of labels of the children of the root is a partition 
of $X$.

\begin{lemma} \label{tb} Let $G$ act on $X$, and let $\tau$ be an
assembly tree for $X$ that is fixed by $G$. If $U$ is the set of
children of the root of $\tau$, then
$U$ is a block system that is compatible with the action of $G$ in $X$. 
\end{lemma}

\proof For any $v\in U$, let $\tau_v$ be the rooted, labeled subtree of
$\tau$ that consists of root $v$ and all its descendents.  If $\tau$ is
fixed by $G$, then $g(\tau) = \tau$ for each $g\in G$.  
In other words,  $g(\tau_v)= \tau_u$ for some $u\in U$.  
This implies that $g(v)\cap v = \emptyset$ if $u\ne v$ or 
$g(v)\cap v = v$ if $u = v$.
Hence $U$ is  block system that is compatible with the action
of $G$ on $X$.  \qed \B

The following notation will be used in this section. The set of orbits
of $G$ acting on $X$ will be denoted by ${\bf O}$.  For $H \leq G$,
let ${\bf C}_H$ denote a set of (say left) coset representatives of $H$ 
in $G$.
Note that $|{\bf C}_H| = (G:H)$.  For $r \in G$ and $Q
\subseteq X$, let $r(Q) := \{ r(q) \, : \, q\in Q\}$.
A group $G$ is said to act {\it simply} on $X$ if the stabilizer of 
each $x$ in $X$ is the trivial group.
In this paper, a {\em partition} $\Pi$ of a finite set $S$ into $k$ parts
is a set  $\{\pi_1,\pi_2,\cdots, \pi_k\}$
of disjoint subsets 
so that $\cup_{i=1}^k \pi_i =S$. The subsets $\pi_i$ are called the
{\em parts} of the partition $\Pi$. 
The order of the parts of a partition is insignificant. That is, 
$\{ \{1,3\},\{2,4\} \}$ and $\{\{2,4\}, \{(1,3\} \}$ are identical 
partitions of the set $\{1,2,3,4\}$. We nevertheless label the parts
from $1$ to $k$ for convenience. 

\begin{theorem} \label{CC} Let us assume that $G$ acts simply on $X$.
Let $\Pi = \{\pi_1, \pi_2, \dots, \pi_k\}$ be a partition of ${\bf O}$
into arbitrarily many parts,
 and let ${\bf H} = \{H_1, H_2, \dots,
H_k\}$ be a corresponding sets  of subgroups of $G$.  
For each $i$ and each
$O\in \pi_i$, let $Q_{i,O}$ be any single orbit of the simple
 action of $H_i$ on $O$. 
Let $Q_i = \cup_{O\in \pi_i} Q_{i,O}$ and ${\bf Q} = \{Q_1,Q_2,
\dots , Q_k\}$. Let us denote by $(\Pi,{\bf H},{\bf Q})$ the arrangement 
$\{(\pi_1,H_1,Q_1),\ldots,(\pi_k,H_k,Q_k)\}$
of each $\pi_i$ in $\Pi$ with 
a corresponding subgroup $H_i \le G$ in ${\bf H}$, and $Q_i\in {\bf Q}$.

\begin{enumerate} 
\item The collection $${\bf B}(\Pi,{\bf H},{\bf Q}) = 
\bigcup\limits_{i=1}^k \,
\bigcup\limits_{r \in {\bf C}_{H_i}} 
 r(Q_i),$$ of blocks $r(Q_i)$
is a compatible block system for $G$ acting on $X$.
\item Every compatible block system for $G$ acting on $X$ is
of the above form for some choice of $\Pi, \,{\bf H}$, and $\bf Q$.
\item Two such block systems ${\bf B}(\Pi, {\bf H},{\bf Q})$ and
${\bf B}(\Pi', {\bf H}',{\bf Q'})$ are equal if and only if, there is
a permutation $p$ of the set of blocks of $\Pi'$ so that for all 
$i\leq k$, we have 
$\pi_{p(i')} = \pi_i$, and for all $i\leq k$,
 there exists a $g_i\in G$ so that 
$H_{p(i')} = g_iH_ig_i^{-1}$, and $Q_{p(i')} = g_i(Q_i)$.
\end{enumerate}
\end{theorem}

{\it Proof:} In order to prove Statement (1), let $H \in
{\bf H}$ and $B = r(Q)$, where $Q=Q_i$ for some $i$. We first show that $B$
 is a block.  There is a subset $A \subseteq X$ containing at most one 
element from each $G$-orbit such that $B=rH(A)$.  If 
$g(B) \cap B \neq \emptyset$, then
there are elements $a, a'\in A$ such that $grh(a) = rh'(a')$ for some
$g\in G$ and $h,h'\in H$.  Thus $a$ and $a'$ are in the same
$G$-orbit, which implies that $a=a'$.  Therefore, $grh(a)=rh'(a)$.
Since $G$ acts simply, this implies that  $grh=rh'$, which in turn
implies that $gr$ and $r$ are in the same coset of $H$ in $G$. Therefore
$g(B) = gr(Q) =r(Q) = B$.  This proves,
not only that $B$ is a block, but that ${\bf B}(\Pi,{\bf H},{\bf
Q})$ is a block system, because ${\bf B}(\Pi,{\bf H},{\bf
Q})$ is a partition of $X$ into blocks. Moreover, if $r(Q) \in 
{\bf B}(\Pi,{\bf H},{\bf Q})$ and $g\in G$, then by
definition $gr(Q) \in {\bf B}(\Pi,{\bf H},{\bf Q})$, which 
shows that  ${\bf B}(\Pi,{\bf H},{\bf Q})$ is a block
system compatible with $G$.
 
In order to prove Statement (2), 
let us denote the set of orbits of $G$ in its action on $X$
 by $\{O_1,O_2,\dots, O_n\}$.  
We first show that any block
$B$ in the action of $G$ on $X$ is of the form $B = B_1 \cup B_2 \cup
\cdots \cup B_n$, where $B_i$ is a single orbit of some subgroup
$H\leq G$ acting on $O_i$. Let $B_i = B \cap O_i$. Note that $B_i =
\emptyset$ is a possibility, in which case we have $B = B_1 \cup B_2
\cup \cdots \cup B_m,\, m\leq n$.  Each $B_i$ itself must be a block
because, if $g(B_i) \cap B_i \neq \emptyset$, then $g(B) \cap B \neq
\emptyset$. However, $B$ is a block, so $g(B) \cap B \neq \emptyset$
implies that $g(B) = B$, and thus $g(B_i) = B_i$.

Let $H_i = \{ h\in G \, | \, h(B_i) = B_i\}$. We claim that $H_1 = H_2
= \dots = H_m$.  To see this let $h\in H_i$. Since $B$ is a block,
either $h(B) \cap B = \emptyset$ or $h(B) = B$.  However, $h(B) \cap B =
\emptyset$ is impossible because $h(B_i) = B_i$.  Hence $h(B) =
B$. Now $B_j = B \cap O_j$ implies, for each $j$, that $h(B_j) = B_j$.
Therefore $h_i \in H_j$ for all $i,j$.  This verifies the claim, so
let $H = H_1 = H_2 = \dots = H_m$. 

The proof that each block $B$ is of the required form is complete if
it can be shown that $H$ acts transitively on $B_i$ for each $i$.  To
see this, let $x,y\in B_i$.  Since $B_i$ lies in a single $G$-orbit,
there is a $g\in G$ such that $g(x) = y$.  Since $B_i$ has been shown
to be a block and $g(B_i) \cap B_i \neq \emptyset$, it must be the
case that $g(B_i) = B_i$.  Therefore $g \in H_i = H$.

To complete the proof of Statement (2), let ${\bf B}$ be any
compatible block system for $G$ acting on $X$.  We have proved that
if $B \in {\bf B}$, then $B = B_1 \cup B_2 \cup \cdots \cup B_m$,
where $B_i$ is a single orbit of some subgroup $H\leq G$ acting on
$O_i$.  Because of the compatibility, the action of $G$ on $X$ induces
an action of $G$ on ${\bf B}$. The orbits under this action provide a
partition $\Pi$ of ${\bf O}$, a part $\pi \in \Pi$ consisting of all
$G$-orbits acting on $X$ contained in the union of a single $G$-orbit
acting on ${\bf B}$.  Consider any orbit $W$ of $B$ in this
action. If $B'$ is another element of $W$, then there is an $r\in G$
such that $B' = r(B)$.  This shows that the blocks in $G(B)$ are of
the desired form in Statement (1) of the theorem.  Repeating this
argument for each part in the partition $\Pi$ completes the proof of
Statement (2).

To prove Statement (3), 
we first show that if $Q\in {\bf Q}$ is the
union of $H$-orbits and $H'(Q) = Q$, where $H, H' \in {\bf H}$, then
$H' = H$.  Restricting attention to just one orbit of $G$ in its
action on $X$, the equality $H'(Q) = Q$ implies that $H'(a') = H(a)$
for some $a,a'$ in the same $G$-orbit acting on $X$.  Let $g\in
G$ be such that $a=g(a')$ and hence $Hg(a') = H'(a')$, which
in turn implies that $hg(a') = a'$ for some $h\in H$. Because $G$
acts simply, this implies that $g = h^{-1} \in H$, so $H(a') = H'(a')$,
which again, by the simplicity of the action, implies that $H'=H$.

Now let us assume that ${\bf B}(\Pi, {\bf H} ,{\bf Q})
= {\bf B}(\Pi',{\bf H}',{\bf Q'})$.  Clearly, $\Pi = \Pi'$.
  It is sufficient to restrict our attention to just
one of the parts in the partition $\Pi = \Pi'$, so we must show that
$\{ r(Q) \, : \, r \in {\bf C}_{H} \} = \{ r(Q') \, : \, r \in {\bf
C}_{H'} \}$ if and only if $H' = gHg^{-1}$, and $Q' = g(Q_i)$ for some
$g \in G$.  If $H' = gHg^{-1}$, and $Q' = g(Q)$ for some $g \in G$,
then for any $r \in G$ we have $r( Q') = rH'(Q') =(rgHg^{-1}) g(Q) =
rgH(Q)$.  This shows that $\{ r(Q') \, : \, r \in {\bf C}_{H'} \}
\subseteq \{ r(Q) \, : \, r \in {\bf C}_{H} \}$, and the opposite
inclusion is similarly shown. Conversely, assume that $\{ r(Q) \, : \,
r \in {\bf C}_{H} \} = \{ r(Q') \, : \, r \in {\bf C}_{H'} \}$.
Since $Q' \in \{ r(Q') \, : \, r \in {\bf C}_{H'} \}$, we know that
$Q' = r(Q)$ for some $r \in {\bf C}_{H} \subseteq G$.  
Now $(rHr^{-1})(Q') = (rHr^{-1})(r(Q)) = rH(Q) = r(Q) = Q'$. By the
uniqueness result shown in the preceding paragraph, we get 
  $H' = rHr^{-1}$.  \qed \m

\begin{example} \label{exBlocks} 
Klein 4-group acting on ${\mathcal T}_4$ (continued). 
\end{example}

Continuing the example from the previous section with $G = \Z _2
\oplus \Z _2$ acting simply on $X = \{1,2,3,4\}$, let
$K = K_1 = \{ (1)(2)(3)(4), (1\, 2)(3 \, 4) \}$
and  $K_0$ the trivial subgroup.  
There are 11 blocks in the action of $K$ on $X$ which are given below: 
$$\{1,2,3,4\}, \{1,2\}, \{3,4\}, \{1,3\}, \{2,4\}, \{1,4\}, \{2,3\},
\{1\},\{2\}, \{3\},  \{4\}.
$$ The seven block systems for the action of $K$ on $X$ can be found using
Theorem~\ref{CC}.  In what follows, $\{1,2\} | \{3,4\}$ denotes the orbit
$\{1,2\},\{3,4\}$ partitioned into the two parts $\{1,2\}$ and $\{3,4\}$, 
whereas $\{1,2\},\{3,4\}$ denotes that same orbit partitioned the trivial way,
into one part. Note that ${\bf B}\;(\;\{1,2\} | \{3,4\} ,\;
\{K_0,K\}\;,\; \{2\}|\{3,4\}\;)$, for example, is not included in
the list below.  This is because, according to Statement (3) in
Theorem~\ref{CC}, $$\begin{aligned} 
& {\bf B}\;(\;\{1,2\}|\{3,4\} ,\; \{K_0,K\}\;,\;
\{2\}|\{3,4\}\;) = \\
& {\bf B}\;(\;\{1,2\}|\{3,4\} ,\;
\{K_0,K\}\;,\; \{1\}|\{3,4\}\;).\end{aligned}$$ 
Namely, for $g = (1\,2)(3\, 4)$,
we have $\{2\} = g(\{1\})$ and $K_0 = gK_0g^{-1}$. 

$$ \begin{aligned} 
{\bf B}\;(\;\{1,2\},\{3,4\} ,\; \{K\},\; \{1,2,3,4\}\; \;
& = \; (1 \, 2 \, 3 \, 4) \\
 {\bf B}\;(\;\{1,2\},\{3,4\} , \;\{K_0\},\; \{1,3\}\;) & = \; (1 \,
3)(2 \, 4) \\ 
 {\bf B}\;(\{1,2\},\{3,4\} ,\; \{K_0\},\; \{1,4\}\;) & = \; (1 \,
4)(2 \, 3) \\ 
 {\bf B}\;(\;\{\{1,2\}|\{3,4\} ,\; \{K,K\}, \;\{1,2\},\{3,4\}\;) & =
 \; (1 \, 2)(3 \, 4) \\
{\bf B}\;(\{\{1,2\}|\{3,4\} ,\; \{K,K_0\},\; \{1,2\},\{3\}\;) & =
 \; (1 \, 2)(3)(4) \\
{\bf B}\;(\{1,2\}|\{3,4\} ,\; \{K_0,K\},\; \{1\}, \{3,4\}\;) & =
 \;  (1) (2)(3 \, 4) \\
{\bf B}\;(\{1,2\}|\{3,4\} , \;\{K_0,K_0\},\; \{1\},\{2\}\;) & =
 \;  (1) (2)(3) (4) 
\end{aligned}$$

Let $\tau \in {\mathcal T}_X$ be a tree fixed by $G$ in its action on 
$T_X$.
If $U$ denotes the set of children of the root of $\tau$, 
recall that Lemma~\ref{tb} states that the set $U$ of labels 
is a block system.  Recall that the label of a vertex is the set of labels
of its leaf descendents, and also the label of a vertex is the union of the 
labels of its children. According to Theorem~\ref{CC}, any block system is
of the form
$${\bf B}(\Pi,{\bf H},{\bf Q}) = \bigcup\limits_{i=1}^k \,
\bigcup\limits_{r \in {\bf C}_{H_i}}r(Q_i).$$ 
We will use the notation $\tau_{rQ}$ to
denote the subtree $\tau_u, \, u\in U, u = r(Q)$, rooted at $u$. 
Theorem~\ref{CC} leads to the characterization
of assembly trees fixed by given group $G$ as stated in
Corollary~\ref{a-tree} below.  

\begin{cor} \label{a-tree}
Let us assume
 that $G$ acts simply on $X$ and that $\tau \in {\mathcal T}_X$. Let
$U$ be the set of children of the root of $\tau$ and, for each $u\in U$,
let $\tau_u$ be the rooted, labeled subtree of $\tau$ that consists of root
$u$ and all its descendents. With notation as in Theorem~\ref{CC}, the
tree $\tau$ is fixed by $G$ if and only if, for some $\Pi, \, {\bf H}$ and
$\bf Q$, the following two conditions hold.
\begin{enumerate}
\item $U = {\bf B}(\Pi,{\bf H},{\bf Q})$, hence for each $Q\in {\bf Q}$
and $g\in G$, there is a subtree $\tau_Q$ and a subtree $\tau_{gQ}.$
\item $\tau_{gQ} = g(\tau_{Q})$ for every $Q \in {\bf Q}$
and every $g \in G$.
\end{enumerate}
\end{cor}

\proof Let us assume that $\tau$ is fixed by $G$. Condition (1) follows
immediately from Lemma~\ref{tb} and Theorem~\ref{CC}.  Concerning
Condition (2), for any $g\in G$, the set of leaves of $g(\tau_Q)$ is
$g(Q)$.  Hence for $\tau$ to be fixed by $G$ it is necessary that $g(\tau_Q)
= \tau_{gQ}$.

Conversely, let us assume that Conditions (1) and (2) hold.  For any $g\in G$
we must show that $g(\tau) = \tau$.  By Condition (1), it is
sufficient to show that $g$ acting on $\tau$ permutes the set of subtrees
in such a manner that $g(\tau_{rQ}) = \tau_{grQ}$ for every $H \in
{\bf H}$, $\; Q$ the corresponding element of $\bf Q$,  and every $r
\in {\bf C}_{H}$.  However, by Condition (2), 
$g(\tau_{rQ}) = gr(\tau_Q) = \tau_{grQ}$. 
\qed \B

Theorem~\ref{alg} below states that the following recursive
procedure constructs any assembly tree $\tau \in {\mathcal T}_X$ fixed
by $G$.  This will be used to prove Theorem~\ref{generating} in
the next section.

\begin{alg} \label{const}
Recursive construction of any assembly tree fixed by a group $\mb G$:
\end{alg}

\begin{enumerate}
\item[(1)] Partition the set ${\bf O}$ of $G$-orbits of $X$: 
$\Pi = \{\pi_1, \pi_2, \dots, \pi_k\}$. Note that the parts of $\Pi$
are labeled $1,2,\cdots ,k$ in some arbitrary way.
\item[(2)] For each $i=1,2,\dots, k$, choose a subgroup $H_i \leq G$.
(If $\Pi$ has only one part then $H_i=G$ is not allowed.)
\item[(3)] For each $i$, choose a single orbit of $H_i$ acting on each of the 
$G$-orbits in $\pi_i$, and let $Q_i$ be the union of these $H_i$-orbits.
\item[(4)] Recursively, let $\tau_{Q_i}$ be any rooted tree whose 
leaves are labeled by $Q_i$ 
and which is fixed by $H_i$.
\item[(5)] 
Let  $S_i = \{\ r(\tau_{Q_i})\  |\   r \in {\bf C}_{H_i}\}$
and $S = \cup_{i=1}^k \, S_i$.  Let $\tau$ be the rooted tree whose children
 are roots of the trees in $S$.
\end{enumerate}

\begin{theorem} \label{alg} The set of assembly trees constructed by 
Procedure~\ref{const} is the set of assembly trees fixed by the group
$G$. 
\end{theorem}

\proof In the notation of Theorem~\ref{CC}, Steps (1), (2), and (3) are
choosing $(\Pi,{\bf H},{\bf Q})$. Steps (4), (5), and (6)
are ensuring that $U = {\bf B}(\Pi,{\bf
H},{\bf Q})$.  
Note that the restriction in Step (2) is because otherwise the
root of the resulting tree in Step (6) would have only one child.
Note also that in Step (5), $S_i$ does not depend on the particular
set of coset representatives.  This follows directly from Step (4).

It is now sufficient to show the following. For any assembly tree $\tau$
satisfying $U = {\bf B}(\Pi,{\bf H},{\bf Q})$
for some $(\Pi,{\bf H},{\bf Q})$,
Condition (2) in Corollary~\ref{a-tree} holds if and only if $\tau$ is
constructed by Procedure~\ref{const}.  To show that any assembly tree
$\tau$ constructed by Procedure~\ref{const} satisfies Condition (2), note
that Step (4) implies that, if $Q \in {\bf Q}$ corresponds to $H\in
{\bf H}$, then $H(Q) = Q$ and hence $h(\tau_Q) = \tau_{hQ} = \tau_Q$ for all 
$h\in H$.  For $g\in
G$, if $g=rh$, where $h\in H$, then $g(\tau_Q) = rh(\tau_Q) =
r(\tau_{hQ}) = r(\tau_Q) = \tau_{rQ}$, the last equality from Step (5).
Again, because $H(Q) = Q$, we have $g(\tau_Q) = \tau_{rQ} = \tau_{rhQ} =
\tau_{gQ}$.

Conversely, if $\tau$ satisfies Condition (2) in Corollary~\ref{a-tree},
then consider the trees $\tau_{Q_i}$ 
$i =1,2,\dots , k$. These are trees whose leaves are labeled by $Q_i$.
in Step (4) of Procedure~\ref{const}. 
Moreover, by Condition (2) we have $h(\tau_{Q_i}) 
=  \tau_{hQ_i} = \tau_{Q_i}$ for all $h \in H_i$, so $\tau_{Q_i}$ is
fixed by $H_i$.  By Step (5) of Procedure~\ref{const} and Condition (2)
of  Corollary~\ref{a-tree} we have $ r(\tau_{Q_i}) =
\tau_{rQ_i}$ for all $r\in {\bf C}_{H_i}$.  Therefore the tree $\tau$ is
constructed by  Procedure~\ref{const}.  
 \qed \m

\begin{remark} \label{unique} Enforcing uniqueness in
 the construction. \end{remark} The construction in Procedure \ref{const} is 
not unique, in that  it may produce the same fixed
assembly tree multiple times depending on the choices in Steps 2 and 3.  
 Condition (3) in Theorem~\ref{CC} shows that 
we may enforce uniqueness  if we make the following two restrictions.

\begin{enumerate}
\item[(a)] 
If we choose $(\Pi,{\bf H},{\bf Q}) = \{(\pi_1,H_1,Q_1),\ldots,
(\pi_k,H_k,Q_k)\}$
in Steps 1 and 2 of the procedure while constructing a tree $\tau$,
and if we also have 
$(\Pi,{\bf H'},{\bf Q'}) = \{(\pi_1,H_1',Q_1'),\ldots,(\pi_k,H_k',Q_k')\}$
during the construction of another tree $\tau'$, then to ensure that
$\tau \ne \tau'$ we need to ensure that for at least one $i$,
the group $H_i'$ should not be conjugate to $H_i$ in $G$.
\item[(b)] Consider the construction of two trees $\tau$ and $\tau'$ 
with corresponding 
$(\Pi,{\bf H},{\bf Q})$ and $(\Pi,{\bf H'},{\bf Q'})$ such that for each
$i$, the subgroup $H_i$ is a conjugate of the subgroup $H_i'$. Further
 assume that in Step (3) for the tree $\tau$ 
the element $g_i\in H_i$ is such that $g_iH_ig_i^{-1}=H_i'$.
 Then while constructing 
 tree 
$\tau'$, we need to ensure that there is at least one $i$ 
such that $Q_i'\neq g_i(Q_i)$.  (Note that for a given 
index $i$, there may well be several elements $g_i\in G$
 so that $g_iH_ig_i^{-1} = H_i'$
holds, and all those  are subject to  this 
restriction.)

\end{enumerate}

\begin{example} \label{exFixed} 
Klein 4-group acting on ${\mathcal T}_4$ (continued).
\end{example}

With $G = \Z _2 \oplus \Z _2$ acting on $X = \{1,2,3,4\}$, consider the
assembly trees $\tau$ fixed by the subgroup $K = \{ (1)(2)(3)(4),
(1\,2)(3\, 4) \}$.  There are exactly six such trees, those in the
orbits $A,B,C,D,E$ of Figure~\ref{trees4}. These correspond 
(not in corresponding order) to the
block systems in Example~\ref{exBlocks}.  Because of the restriction
in Step (2) of Procedure~\ref{const}, the first block system in the
list in Example~\ref{exBlocks} is ignored.

\subsection{When the action of  $G$ on $X$ is not simple} 

Let us assume that $G$ acts on $X$, but not necessarily simply.  For $q \in
X$, let $S_q$ denote the stabilizer of $q$ in $G$. For a subset $Q
\subseteq X$, let $$S_Q = \bigcup_{q\in Q} S_q.$$ If $G$ acts simply
on $X$, then the stabilizer of any $x\in X$ is the trivial
subgroup. Therefore, in this case, it is clear that $S_Q \subset H$
for any $Q\subseteq X$ and $H\leq G$.  In the general case, when
$G$ acts not necessarily simply on $X$, let us call a pair $(H,Q)$ {\it
viable} if
$${\bf S}_Q \subseteq H.$$ If only viable pairs $(H_i,Q_i)$ are
allowed in the hypothesis of Theorem~\ref{CC}, then the theorem is
valid in the general, not necessarily simple, case.  Since this
general version of Theorem~\ref{CC} and associated analogs
of Procedure~\ref{const} and Theorem~\ref{alg} 
are not needed in subsequent
sections, and the proofs are relatively straightforward extensions, 
we omit them. 

\section{Enumerating Fixed Assembly Trees}
\label{enumerating}

Let us assume in this section that $G$ acts simply on each of an infinite
sequence $X_1, X_2, \dots$ of sets where, by formula (\ref{orbitsize})
we have $|X_n| = n|G|$.  In other words, $n$ is the number of orbits
of $G$ in its action on $X_n$.  Denote by $t_n(G)$ the number of trees
in ${\mathcal T}_n := {\mathcal T}_{X_n}$ that are fixed by $G$.  In this
section we provide a formula for the exponential generating function
$$f_G(x) := \sum_{n\geq 1} t_n(G) \,\frac{x^n}{n!}$$ for the sequence
$\{t_n(G)\}$.  If $G$ is the trivial group of order one, then let us
 denote this
generating function simply by $f(x)$. This is the generating function
for the total number of rooted, labeled trees with $n$ leaves in which
every non-leaf vertex has at least two children. For $H\leq G$, let
$$\widehat f_H (x) = 
\frac{1}{(G:H)}\; f_H\left( (G:H)x \right).$$

\begin{theorem} \label{generating} The generating function
$f_G(x)$ satisfies the following functional equations:
$$1-x+2f(x) = \exp \, (f(x)),$$ and for $|G|>1$, $$1+2 f_G(x) =
\exp\, \left(\sum_{H \leq G} \; \widehat f_H (x)\right).$$
\end{theorem}

\proof The first formula is proved in \cite{Stanley}, page 13.  For
$|G|>1$, we use the standard exponential and the product formulas
for generating functions.

The proof of the second formula uses two well known results from the
theory of exponential generating functions, the ``product formula''
and the ``exponential formula''. In Procedure~\ref{const}, give 
Steps (3) and (4) the name {\it putting an $H_i$-structure on
$\pi_i$}. According to Theorem~\ref{alg}, the number of trees $t_n(G)$
fixed by $G$ equals the number of ways to partition the set of orbits
of $G$ acting on $X_n$ and to place an $H$-structure on each part in
the partition, for some subgroup $H \leq G$, keeping the uniqueness
Remark~\ref{unique} in mind.   

In Step (3) of Procedure~\ref{const}, since $G$ acts simply and the
number of $H_i$-orbits in one $G$-orbit is $|G|/|H_i| = (G:H_i)$, the
number of possible choices for $Q_i$ (the union of these single
$H_i$-orbits) is $(G:H)^m$.  Hence, in accordance with Step (4) of
Procedure~\ref{const}, the generating function for the number of ways
to place an $H$-structure is basically $f_H((G:H)x)$.  

However, this must be altered in accordance with the uniqueness
requirements in Remark~\ref{unique}.  Let ${\bf N}$ denote a set
consisting of one representative of each conjugacy class in the set of
subgroups of $G$.  By Statement (a) in Remark~\ref{unique}, only
subgroups in ${\bf N}$ are considered. Let $N(H) := \{g\in G \, | \,
gHg^{-1} = H \}$ denote the {\it normalizer} of $H$ in $G$. 
 By
Statement (b),  there has to be an index $i$ so that $g(Q_i) \neq Q_i'$.
However, $g(Q_i)=g'(Q_i)$ will occur for every $i$ if and only if
$g$ and $g'$ are in the same coset of $H$ in $G$.
 Therefore, the generating function for the
number of ways to place an $H$-structure is $\frac{1}{(N(H):H)} \;
f_H((G:H)x)$.

The exponential formula states that the generating function $g_H(x) =
\sum_{n\geq 0} a_n \,\frac{x^n}{n!}$ for the number of ways $a_n$ to
partition the set of $G$-orbits acting on $X_n$ and, on each part
$\pi$ in the partition, place an $H$-structure (same $H$) is
$$g_H( x) := \exp \left( \widehat f_H (x)\right ).$$ Here we assume 
that $a_0 = 1$.  
 
The generating function for the number of ways to partition the set of
orbits, i.e., choose $\Pi = (\pi_1, \pi_2, \dots, \pi_k)$ and, on each
part of the partition, place an $H$-structure, one $H$ from each
conjugacy class in ${\bf N}$ is
$$\begin{aligned}
\prod_{H \in {\bf N}} g_H(x) &= \prod_{H \in {\bf N}} 
\exp \left(
\frac{1}{(N(H):H)} \; f_H((G:H)x) \right )\\ 
&= \exp\, \left(\sum_{H \in {\bf N}} \; \frac{1}{(N(H):H)} \;
f_H((G:H)x) \right). \end{aligned}$$
Note that we have not taken the restriction in Step (2) of
Procedure~\ref{const} into consideration.  Taking the partition of the
orbit set into just one part and placing on that part a $G$-structure
results in counting the number of fixed trees a second time. Also  
since the constant term in $\prod_{H \in {\bf N}} g_H(x)$ is 1,
$$\begin{aligned}
1 + 2\,f_G(x) &= \exp\, \left(\sum_{H \in {\bf N}} \; \frac{1}{(N(H):H)} \;
f_H((G:H)x) \right) \\ & = 
\exp\, \left(\sum_{H \leq G} \; \frac{1}{(G:H)} \;
f_H((G:H)x) \right). \end{aligned}$$
Here the last equality holds because $f_H(x)$ depends only the conjugacy class
of $H$ in $G$ and  
$$ \frac{1}{(N(H):H)} \Big/ \frac{1}{(G:H)} = \frac{(G:H)}{(N(H):H)}
= (G:N(H)) = |{\bf N}|.$$ \qed 

\begin{example}  Klein 4-group acting on ${\mathcal T}_4$ (continued).
\end{example}

Consider $G = \Z_2 \oplus \Z_2$ acting on $X_n$. Recall that $X_n =
4n$, the integer $n$ being the number of $G$-orbits.  In this case
${\bf N} = \{ K_0, K_1,K_2,K_3, G \}$, where $K_0$ is the trivial
group and
$$\begin{aligned}
K_1 &= \{ \, (1)(2)(3)(4), (1\;2)(3\;4) \, \}, \\
K_2 &= \{ \, (1)(2)(3)(4), (1\;3)(2\;4) \, \}, \\
K_3 &= \{ \, (1)(2)(3)(4), (1\;4)(2\;3) \, \}.
\end{aligned}$$
The functional equations in the Statement of Theorem~\ref{generating}
are
$$\begin{aligned}
1-x+2f(x) &= \,\exp \, (f(x)) \\
1+2f_{K_i}(x) &=\, \exp \, \left(\frac12 \,f(2x)+ f_{K_i}(x) \right) 
\quad  \text{for $i=1,2,3$, and}\\
1+2f_G (x) &= \,\exp \,  \left(\frac14 \, f(4x) +\frac12 \, f_{K_1}(2x)+
\frac12 \, f_{K_2}(2x)+\frac12 \, f_{K_3}(2x)+ f_G (x) \right). 
\end{aligned}$$
Using these equations and MAPLE software, the coefficients of the
respective generating functions provide the following first few values
for the number of fixed assembly trees. For the first entry $t_1(G)
=4$ for the group $G$, the four fixed trees are shown in
Figure~\ref{trees4} $A$, $B$, $C$, $D$.  For trees with eight leaves there are
$t_2(G) =104$ assembly trees fixed by $G = \Z_2 \oplus \Z_2$, and so
on.
$$\begin{aligned}
t_n(K_0) \quad &: \quad  1, 1 , 4, 26, 236, 2752  \\
t_n(K_i) \quad &: \quad  1, 6, 72, 1312,  32128, 989696 \\
t_n(G) \quad &: \quad 4, 104, 4896, 341120, 31945728, 3790876672.
\end{aligned}$$

Theorem~\ref{generating} provides the generating function for the
numbers $t_n(H)$ of fixed assembly trees in the action of any subgroup
$H\leq G$ on $X_n$.  What is required for Problem (i) described in
Sections 1 and 2 are the numbers $\overline t_n(H)$ of assembly
trees that are fixed by $H$, but by no other elements of $G$.  In
Example~\ref{exFixed}, for $G=\Z_2 \oplus Z_2$ acting on $X =
\{1,2,3,4\}$, there are six trees that are fixed by the subgroup $K =
\{\, (1)(2)(3)(4) , (1\;2)(3\;4) \, \}$. However, of these six, four 
 ($A$, $B$, $C$, and $D$ in Figure~\ref{trees4}) are also fixed by $G$.
Therefore there are only two assembly trees fixed by $K$ and no other
elements of $G$ (these are $E$, and $F$
 in Figure~\ref{trees4}).  In general, as shown Theorem \ref{mobius},
M\"obius
inversion \cite{VW} can be used to calculate the values of $\overline
t_X(H)$ from the values of $t_X(H)$.

\section{The Icosahedral Group} \label{icos}

For completeness, the results of the previous sections are applied to
the motivating $T=1$ viral example.
 An {\it isometry} of 3-space is a bijective transformation that
preserves length, and an isometry is called {\it direct} if it is
orientation preserving.  Rotations, for example, are direct, while
reflections are not.  A {\it symmetry} of a polyhedron is an isometry
that keeps the polyhedron, as a whole, fixed, and a {\it direct
symmetry} is similarly defined.  The {\it
icosahedral group} is the group of direct symmetries
of the icosahedron. It is a group of order 60 denoted $G_{60}$.

As mentioned earlier, 
the viral capsid is modeled by a polyhedron $P$ with icosahedral symmetry,
whose set $X$ of
facets represent the protein monomers. The icosahedral group, 
 acts on $P$ and hence on the set $X$.  It follows
from the quasi-equivalence theory of the capsid structure that $G_{60}$ acts
simply on $X$.  Formula (\ref{orbitsize}) shows  that 
$|X| := |X_n| = 60n$, where
$n$ is the number of orbits.  Not every $n$ is possible for a viral
capsid; $n$ must be a $T$-number as defined in the introduction.
Before the number of orbits of each size for the action of $G_{60}$ on
the set ${\mathcal T}_n:= {\mathcal T}_{X_n}$ of assembly trees can be
determined, basic information about the icosahedral group is needed. \m

The group $G_{60}$ consists of:

\begin{itemize}
\item the identity,
\item 15 rotations of order 2 about axes that pass through the
midpoints of pairs of diametrically opposite edges of $P$,
\item 20
rotations of order 3 about axes that pass through the centers of
diametrically opposite triangular faces, and
\item 24
rotations of order 5 about axes that pass through diametrically
opposite vertices.  
\end{itemize}

There are 59 subgroups of $G_{60}$ that play a
crucial role in the theory.  Besides the two trivial subgroups, they
are the following:

\begin{itemize}
\item 15 subgroups of order 2, each generated by one of the
rotations of order 2,
\item 10 subgroups of order 3, each generated by one of the
rotations of order 3,
\item 5 subgroups of order 4, each generated by rotations of order
2 about perpendicular axes,
\item 6 subgroups of order 5, each generated by one of the
rotations of order 5,
\item 10 subgroups of order 6, each generated by a rotation of
order 3 about an axis L and a rotation of order 2 that reverses
L,
\item 6 subgroups of order 10, each generated by a rotation of
order 5 about an axis L and a rotation of order 2 that reverses
L,
\item 5 subgroups of order 12, each the symmetry group of a
regular tetrahedron inscribed in $P$.
\end{itemize}

\begin{figure}[htb]
\vskip 2mm
\begin{center}
\includegraphics[width=3in]{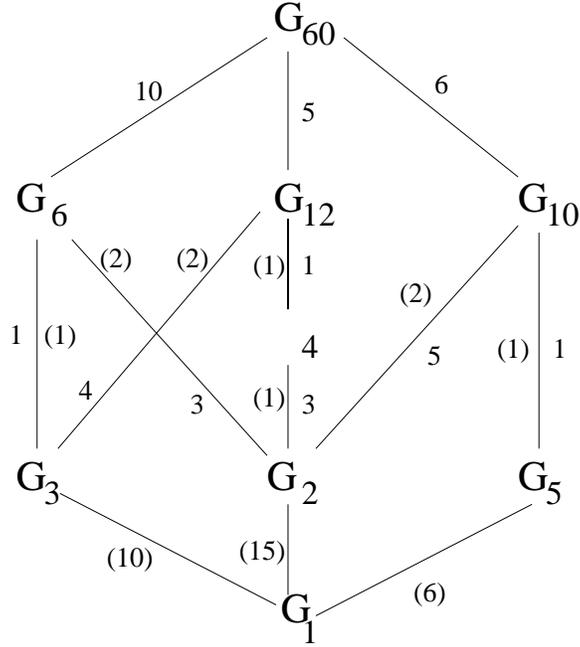}
\end{center}
\caption{Partial Hasse diagram for the lattice of subgroups of the
icosahedral group.} 
\label {Hasse}
\end{figure}

\noindent From the above geometric description of the subgroups, it
follows that all subgroups of a given order are conjugate in the group
$G_{60}$.  Representatives of the conjugacy classes of the subgroups
of the icosahedral group are denoted by $G_0, G_2, G_3, G_5, G_6,
G_{10}, G_{12}, G_{60}$, where the subscript is the order of the
group. The set of subgroups of $G_{60}$ forms a lattice, ordered by
inclusion.  A partial Hasse diagram for this lattice $\bf L$ is shown
in Figure~\ref{Hasse}.  The number on the edge joining $G_i$ (below) and
$G_j$ (above) indicate the number of distinct subgroups of order $i$
contained in each subgroup of order $j$.  The number in parentheses on
the edge joining $G_i$ (below) and $G_j$ (above) indicate the number
of distinct subgroups of order $j$ containing each subgroup of order
$i$.  It is well-known that any finite partially ordered set $P$
 admits a 
M\"obius function $\mu \,: \, P \times P \rightarrow \mathbb
Z$.  The M\"obius function of $\bf L$ is shown in Table 1. The
entry in the table corresponding to the row labeled $G_i$ and column
$G_j$ is $\mu(G_i,G_j)$.

\begin{table}[htb] 
 \begin{center}
  \epsfig{file=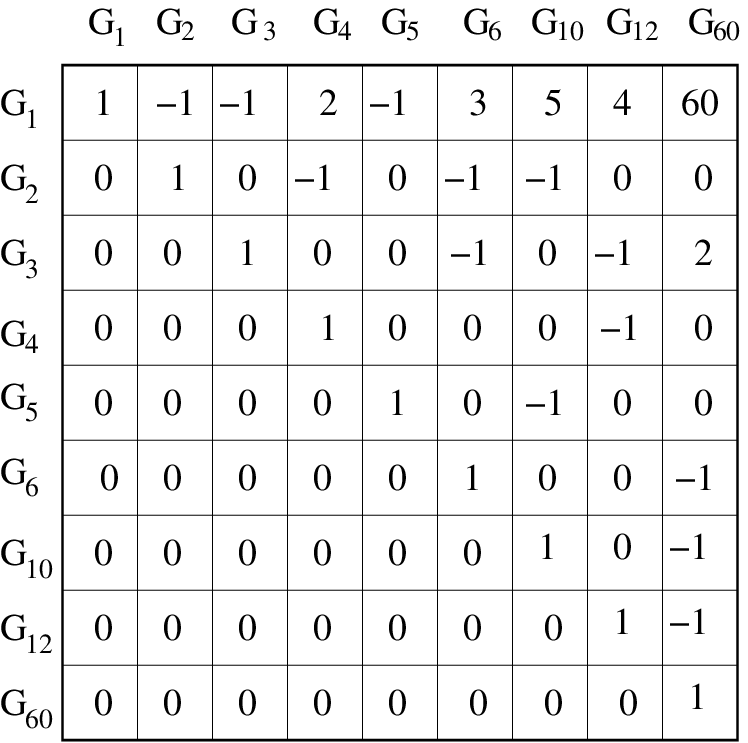}
\caption{The values of the M\"obius function of the subgroup lattice
of $G_{60}$.}
 \end{center}
 \label{mob}
\end{table}

For $|X| = 60$, i.e., for the $T=1$ polyhedral case, 
using Theorem~\ref{generating} and MAPLE software, the generating functions
$f_{G_{i}}(x)$ were computed, and hence their coefficients 
$t_{60/i}(G_{i})$ which
count the number of assembly trees that are 
fixed by $G_{i}$ were also be computed.
Note that since $|X|=60$, the number of orbits of $G_i$ in its action on $X$
is $60/i$. Substituting these values into Theorem~\ref{mobius} and using the 
M\" obius Table 1 yields
the following numerical values for $\overline{t}_i(G_{60/i})$, the 
number of assembly trees over $X$ with $|X| = 60$ that are fixed by $G_i$
but by no other elements of $G_{60}$. In other words, these
are the numbers of  trees whose stabilizer in $G_{60}$
is $G_i$.
$$\begin{aligned}
{\overline t}_{60}(G_{1}) &= 
1924465510132437394720184730922187571120346754532 \\
& \quad 2366329965115755432139023628289410324670840066578537680 \\
{\overline t}_{30}(G_{2}) &=
 1670856367100496379411587456529324583988755126499875584 \\
{\overline t}_{20}(G_{3}) &= 10087157294451731428720995944759704 \\
{\overline t}_{15}(G_{4}) &= 10041342673530270014535171213312 \\
{\overline t}_{12}(G_{5}) &= 20540071766413107840 \\
{\overline t}_{10}(G_{6}) &= 61346927354448105268 \\
{\overline t}_6(G_{10}) &= 223503950260 \\
{\overline t}_5(G_{12}) &= 16865654580 \\
{\overline t}_1(G_{60}) &= 204
\end{aligned}$$

From Theorem~\ref{orbit}, the above numbers $\overline{t}_i(G_{60/i})$
tell us the number of assembly trees with orbit size $i$,
 or in other words, trees in
an assembly pathway of size $i$. That is, the 
probability of such a pathway is $i/|\mathcal{T}_X|$.

It is worth comparing the first and last elements of this list. While the
individual pathways belonging to $G_1$ are only 60 times more probable
then those that belong to $G_{60}$, there are about $10^{99}$
 times more of them.


\section{Conclusion and Open Problems} \label{sec:conclusion} 

We have developed an algorithmic and combinatorial approach to a
problem arising in the modeling of viral assembly.  Our
results illustrate, not only that problems arising from structural
biology can be of independent mathematical interest, but also that
mathematical methods have a direct application in
structural biology.

More specifically, we have developed techniques to analyze the
probability of a capsid forming along a given assembly pathway. One
remaining issue is how to extend these techniques to  finding the
probability of {\sl valid} assembly pathways as defined in Section
\ref{intro}. As mentioned earlier, valid assembly trees can be defined
combinatorially, using generalized notions of connectivity of the
polyhedral graph whose facets form the leaves of the tree.  Combining
such graph theoretic restrictions with our techniques will likely
require new ingredients.  A second important issue is how to extend
our techniques to {\sl nucleation} in viral shell assembly.
Mathematically \cite{bonasith}, the problem is to estimate the
proportion of valid assembly trees that have a subtree whose leaves
form a specific subset of facets, for example a trimer or a pentamer,
in the underlying polyhedron.

In addition to the above extensions of the theory, there is scope to
tighten some results of the paper.  For example, a finer complexity
analysis for Algorithm {\bf Stabilizer} could be based on using Sim's
algorithm, strong generating sets, and the
Cayley graph for $G$ as input.

A study of unlabeled trees that are $g$-unfixable may lead to relevant
related results. Let us say that a tree is $g$-{\it unfixable} if there is no
leaf-labeling so that the resulting labeled tree is fixed by the
permutation $g$, and let us say that a tree is
 $G$-{\it unfixable} if it is $g$-unfixable for
every nontrivial element of the group $G$.  These properties are
interesting for at least two reasons. 
First, they clarify the minimum quantifiable
information in a labeled tree that is needed for 
deciding if it is fixed by a group
element $g$: if the underlying unlabeled tree is $g$-unfixable,  
then the information in the labeling is unnecessary to make this decision.
This may lead to efficient algorithms and tight
complexity bounds.  Second, in the language of formal logic, these
properties are likely to be monadic second order expressible
\cite{bib:compton,bib:woods},
permitting the application of limit laws for the asymptotic
probabilities of finite structures satisfying such properties.


\end{document}